\documentclass[a4paper,10pt]{article}
\usepackage[top=3cm,bottom=3cm,left=2.5cm,right=2.5cm]{geometry}
\usepackage{amssymb}
\usepackage{graphicx}
\usepackage{float}
\usepackage{subfigure}
\usepackage{epstopdf}
\usepackage{amsmath,amsthm,amssymb,lineno}
\usepackage{latexsym}
\usepackage{graphicx,booktabs,multirow}
\usepackage{tikz}
\usetikzlibrary{decorations.pathreplacing}
\usetikzlibrary{intersections}
\usepackage{enumerate}
\usepackage{graphicx,booktabs,multirow}
\usepackage{appendix}
\newtheorem{theorem}{Theorem}
\newtheorem{lemma}[theorem]{Lemma}

\begin{document}

\title{Maximum zeroth-order general Randi\'{c} index of orientations of trees, unicyclic and bicyclic graphs with given matching number}
\author{ Jiaxiang Yang, Hanyuan Deng\thanks{Corresponding author: hydeng@hunnu.edu.cn}\\
{\small Key Laboratory of Computing and Stochastic Mathematics (Ministry of Education),}
 \\{\small School of Mathematics and Statistics, Hunan Normal University,}
 \\{\small  Changsha, Hunan 410081, P. R. China.}}

\date{}
\maketitle

\begin{abstract}
The  zeroth-order general Randi\'{c} index $R^{0}_{a}$ of a  digraph $D$  is  the sum of  $(d^{+}_{v})^{a}+(d^{-}_{w})^{a}$ over all arcs $vw$ of $D$,
where  $a$, $d^{+}_{v}$ and $d^{-}_{w}$ are an arbitrary real number, the out-degree  of the vertex $v$ and the  in-degree of the vertex $w$, respectively.
We determine maximum  zeroth-order general Randi\'{c} index of  oriented trees, unicyclic and  bicyclic graphs in terms of matching number and order in this paper.
\end{abstract}

\maketitle

\makeatletter
\renewcommand\@makefnmark%
{\mbox{\textsuperscript{\normalfont\@thefnmark)}}}
\makeatother

\baselineskip=0.25in

\section{Introduction}

Let $G=(V, E)$ be a simple connected graph with the vertex set $V(G)$  and the edge set $E(G)$. We denote by $d_{G}(v)$ ($d_v$ for short) the degree of vertex $v$ in $G$, and $N_{G}(v)$  the neighbors of vertex $v$ in $G$.
  We denote by $G-uv$ (resp. $G+uv$)  the subgraph of $G$  obtained  by deleting (resp. adding) $uv$ with $uv\in E(G)$ (resp. $uv\notin E(G)$), and $G-u$ (resp.  $G+u$)  the subgraph of $G$  obtained by deleting (resp. adding)  vertex $u$  and the edges incident with $u$, where $u\in V(G)$ (resp. $u\notin V(G)$).
 A matching $M$ of the graph $G$ is a subset of $E(G)$  in which  no two edges  share a common vertex. For any  matching $M_1$ of $G$, if  $|M_1| \leq |M|$, then the matching $M$  is  maximum in $G$. The matching number of $G$ is equal to the number of edges
of a maximum matching in $G$. If  vertex $v\in V(G)$
is incident with an edge of $M$, then $v$ is  $M$-saturated. If  $v$ is $M$-saturated for any $v\in V(G)$, then $M$ is a perfect matching in $G$.

Let  $D=(V,A)$ be a digraph with  arc set $A(D)$  and vertex set $V(D)$, where $vu\in A(D)$  is an arc that from vertex $v$  to vertex $u$ in the graph $D$ (  $D$ has no loops).    $d^{+}_u$ (resp. $d^{-}_u$) denotes the out-degree (resp. in-degree) of a vertex $u$,  $d^{+}_{u}=|\{v| uv\in A(D),v\in V(D)\}|$ (resp. $d^{-}_{u}=|\{v| vu\in A(D),v\in V(D)\}|$). If  $u\in V(D)$ and  $d^{+}_u=d^{-}_u=0$, then $u$ is  an isolated
vertex.
If  $d^{+}_u=0$ (resp. $d^{-}_u=0$) for $u\in V(D)$, then $u$ is  a sink vertex (resp. source vertex).
We can replace each edge $uv$ of $G$ by an arc $uv$ or $vu$ and get an oriented graph $D$ which is also called an orientation of $G$. Let $\mathcal{O}(G)=\{D| D$ is a orientation of $G\}$.
$D\in \mathcal{O}(G)$, if any $u\in V(D)$, we have $d^{+}_u=0$ or $d^{-}_u=0$, then $D$ is  a sink-source orientation of $G$. Let $\mathcal{O'}(G)=\{D|D$ is a sink-source orientation of $G\}$.

Topological index has a great influence   in QSPR
and QSAR \cite{1,2}.
 Li and Zheng \cite{3} proposed
 zeroth-order general Randi\'{c} index  defined as
$R^{0}_{a}(G)=\sum\limits_{u\in V(G)}(d_{G}(u))^{a},$
where  $a$ is  an arbitrary real number.
 It is obvious that the zeroth-order general Randi\'{c} index is a VDB topological index.
The problem that what is extremal graphs over the set of important graphs with respect to zeroth-order general Randi\'{c} index  has been studied by researchers in  \cite {4,5,6,7,9,10,11}.
Recently, Monsalve and Rada  proposed  VDB topological indices in digraphs; see \cite {JM211}.

The definition  of zeroth-order general Randi\'{c} index $R^{0}_{a}(D)$ of a digraph is as \cite{JY223}
$$R^{0}_{a}(D)=\frac{1}{2}\sum\limits_{uv\in A}[(d^{+}_u)^{a}+(d^{-}_v)^{a}]=\frac{1}{2}\sum\limits_{u\in V}[(d^{+}_u)^{a+1}+(d^{-}_u)^{a+1}].$$
If $D$ is obtained by   replacing each edge $uv\in E(G)$
with  arcs $uv$ and $vu$, then $d^{+}_D(u)=d^{-}_D(u)=d_G(u)$ and
$$R^{0}_{a}(D)=\frac{1}{2}\sum\limits_{uv\in A}[(d^{+}_D(u))^{a}+(d^{-}_D(v))^{a}]=\frac{1}{2}\sum\limits_{uv\in E}[(d_G(u))^{a}+(d_G(v))^{a}]= \frac{1}{2}\sum\limits_{u\in V}[(d_G(u))^{a+1}]=R^{0}_{a+1}(G).$$

Monsalve and Rada \cite {JM211} solved the problem on extremal values of the Randi\'{c} index over a important digraphs family, such as oriented paths in terms of the order and cycles in terms of the order, digraphs in terms of the order, oriented  hypercubes in terms of the  dimension,  oriented trees in terms of the order,  respectively.
Deng et al. \cite{HD22} solved the problem on extremal values of some VDB topological indices
 over a important digraphs family that is digraphs in terms of the order, such as the harmonic index,the Randi\'{c} index, the Atom-Bond-Connectivity index,the first and second Zagreb indices, the Geometric-Arithmetic index and the sum-connectivity index.
 The problem of finding extremal values  of  VDB topological
indices over a important digraphs family that are oriented trees in terms of the order and oriented graphs  in terms of the order  were solved by  Monsalve and Rada \cite {JM212}.
 The problem of finding extremal values of first Zagreb  index over a important digraphs family that is oriented  unicyclic graphs in terms of matching number  and the order
  was solved by  Yang and  Deng \cite{JY22}.




 In this paper, we discuss orientations of trees, unicyclic and bicyclic graphs in terms of the matching number with maximum value  of  the zeroth-order general Randi\'{c} index.

\section{Preliminary  }

We first establish three useful lemmas.

\begin{lemma}\cite {JM19}\label{lem16}
 Let $G$ be a graph. Then $G$ is a bipartite graph if and only if $|\mathcal{O'}(G)|\geq 1$. Moreover, If $G$ is a connected bipartite graph, then $|\mathcal{O'}(G)|=2$.
\end{lemma}

\begin{lemma}\cite{JY223}\label{lem1}
For any simple connected  graph $G$,
 $D \in \mathcal{O}(G)$, $a\geq 1$, we have
$$R^{0}_{a}(D)\leq \frac{R^{0}_{a+1}(G)}{2}$$  with equality  if and only if $D\in \mathcal{O'}(G)$.
\end{lemma}


\begin{lemma}\label{lem101}
Let $G$ be a graph of order  $n$, $a\geq 1$, $D\in \mathcal{O}(G)$.

(1)\cite{JY223} If $d_{G}(u)=1$ and $u\in V(G)$, $G'=G-u$,  $D'\in \mathcal{O}(G')$ such that $A(D')\bigcap A(D)=A(D')$, $v \in N_{G}(u)$. Then
$$R^{0}_{a}(D)-R^{0}_{a}(D')\leq \frac{1}{2}[1+(d_{G}(v))^{a+1}-(d_{G}(v)-1)^{a+1}]$$
  with equality if and only if   $\max\{d^{+}_{D}(v), d^{-}_{D}(v)\}=d_{G}(v)$.

(2)If $d_{G}(w_1)=1$, $w_2\in N_{G}(w_1)$, $d_{G}(w_2)=2$ and $w_1,w_2\in V(G)$,  $G'=G-\{w_1,w_2\}$,  $D'\in \mathcal{O}(G')$ such that $A(D')\bigcap A(D)=A(D')$, $ N_{G}(w_2)=\{w_1,w_3\}$. Then
$$R^{0}_{a}(D)-R^{0}_{a}(D') \leq\frac{1}{2}[1+2^{a+1}+(d_{G}(w_3))^{a+1}-(d_{G}(w_3)-1)^{a+1}]$$
 with equality if and only if  $d^{+}_{D}(w_3)=d_{G}(w_3), d^{-}_{D}(w_2)=2$; or $d^{-}_{D}(w_3)=d_{G}(w_3), d^{+}_{D}(w_2)=2$.

(3)If $d_{G}(w_1)=2$,$w_3\in N_{G}(w_1)$,  $w_1,w_3\in V(G)$,  $G'=G-w_1w_3$,  $D'\in \mathcal{O}(G')$ such that $A(D')\bigcap A(D)=A(D')$. Then
$$R^{0}_{a}(D)-R^{0}_{a}(D')\leq \frac{1}{2}[2^{a+1}-1+(d_{G}(w_3))^{a+1}-(d_{G}(w_3)-1)^{a+1}]$$
with equality if and only if  $d^{+}_{D}(w_3)=d_{G}(w_3), d^{-}_{D}(w_1)=2$; or $d^{-}_{D}(w_3)=d_{G}(w_3), d^{+}_{D}(w_1)=2$.
\end{lemma}
\begin{proof}
 Let $b_i=d_{D}^{+}(w_i)$, $b_i'=d_{D'}^{+}(w_i)$, $d_i=d_{D}^{-}(w_i)$, $d'_i=d_{D'}^{-}(w_i)$, $i=1,2,3$.
$d_{G}(w_3)=t$.

(2)

  If $w_2w_3\in A(D)$, then $d_3=d_3'+1$, $b_3=b_3'$. Hence,
 \begin{equation*}
\begin{aligned}
R^{0}_{a}(D)-R^{0}_{a}(D')
=&\frac{1}{2}[(b_1)^{a+1}+(d_1)^{a+1}-0+(b_2)^{a+1}+(d_2)^{a+1}-0]\\
&+\frac{1}{2}[(b_3)^{a+1}+(d_3)^{a+1}-(b_3')^{a+1}-(d_3')^{a+1}] \\
\leq &\frac{1}{2}[1+2^{a+1}]+\frac{1}{2}[(b_3)^{a+1}+(d_3)^{a+1}-(b_3)^{a+1}-(d_3-1)^{a+1}] \\
\leq&\frac{1}{2}[1+2^{a+1}+t^{a+1}-(t-1)^{a+1}]
\end{aligned}
\end{equation*}
 with equality if and only if   $d_3=t, b_2=2$.

  Similarly, if $w_3w_2\in A(D)$, then
 \begin{equation*}
\begin{aligned}
R^{0}_{a}(D)-R^{0}_{a}(D')
&  \leq\frac{1}{2}[1+2^{a+1}+t^{a+1}-(t-1)^{a+1}]
\end{aligned}
\end{equation*}
  with equality if and only if  $b_3=t, d_2=2$.
 The lemma holds clearly.

 (3)
  If $w_1w_3\in A(D)$, then $b_1=b_1'+1$, $d_1=d_1'$, $d_3=d_3'+1$, $b_3=b_3'$. Hence,
 \begin{equation*}
\begin{aligned}
R^{0}_{a}(D)-R^{0}_{a}(D')
=&\frac{1}{2}[(b_1)^{a+1}+(d_1)^{a+1}-(b_1')^{a+1}-(d_1')^{a+1}]\\
&+\frac{1}{2}[(b_3)^{a+1}+(d_3)^{a+1}-(b_3')^{a+1}-(d_3')^{a+1}] \\
=&\frac{1}{2}[(b_1)^{a+1}+(d_1)^{a+1}-(b_1-1)^{a+1}-(d_1)^{a+1}]\\
&+\frac{1}{2}[(b_3)^{a+1}+(d_3)^{a+1}-(b_3)^{a+1}-(d_3-1)^{a+1}] \\
 \leq& \frac{1}{2}[2^{a+1}-1+t^{a+1}-(t-1)^{a+1}]
\end{aligned}
\end{equation*}
 with equality if and only if   $d_3=t, b_1=2$.

  Similarly, if $w_3w_1\in A(D)$, then
 \begin{equation*}
\begin{aligned}
R^{0}_{a}(D)-R^{0}_{a}(D')
& \leq \frac{1}{2}[2^{a+1}-1+t^{a+1}-(t-1)^{a+1}]
\end{aligned}
\end{equation*}
  with equality if and only if  $b_3=t, d_1=2$.
 The lemma holds clearly.
\end{proof}

\section{Maximum zeroth-order general Randi\'{c} index of orientations of trees, unicyclic  graphs with given matching number}
\begin{figure}[ht]
\begin{center}
  \includegraphics[width=5cm,height=2cm]{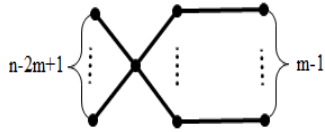}
     \end{center}
\vskip -0.5cm
\caption{ The graph $T_{n,m}$ .}\label{fig-10}
\end{figure}

In this section, we first determine the maximum   zeroth-order general Randi\'{c} index for oriented trees  with given matching number. \\
\indent Let  $1\leq m\leq \lfloor\frac{n}{2} \rfloor$, where $n$, $m$ are integers. Let $T(n,m)=\{T| T$ is a tree  of order $n$ and matching number $m$ $\}$. Let $T_{n,m}$  be shown in Figure \ref{fig-10}.  Obviously,  $T_{n,m}\in T(n,m)$.

Denote $h_0(n,m,a)=\frac{1}{2}[n-m+2^{a+1}(m-1)+(n-m)^{a+1}]$.

\begin{lemma}\cite{LQ10}\label{lem3}
Let $T\in T(n,m)$ with $1\leq m\leq \lfloor\frac{n}{2} \rfloor$, $a\geq 1$. Then
$$R^{0}_{a+1}(T)\leq 2h_0(n,m,a)$$
with equality if and only $T \cong T_{n,m}$.
\end{lemma}

We can extend  the zeroth-order general Randi\'{c} index of trees to the oriented trees.
\begin{theorem}\label{lem4}
Let $T\in T(n,m)$ with $1\leq m\leq \lfloor\frac{n}{2} \rfloor$, $a\geq 1$, $D\in \mathcal{O}(T)$. Then
$$R^{0}_{a}(D)\leq h_0(n,m,a)$$
with equality if and only $D\in \mathcal{O'}(T_{n,m})$.
\end{theorem}
\begin{proof}
     \indent From Lemma $\ref{lem1}$ and Lemma $\ref{lem16}$, $R^{0}_{a}(D)\leq \frac{1}{2}R^{0}_{a+1}(T)$ with equality if and only if $D\in \mathcal{O'}(T)$.  Then by Lemma $\ref{lem3}$, we have $$\max \{R^{0}_{a}(D)|\\D\in \mathcal{O}(T),T\in T(n,m)\}=\max \{\frac{1}{2}R^{0}_{a+1}(T)|T\in T(n,m)\}=\frac{1}{2}R^{0}_{a+1}(T_{n,m}).$$
     Consequently,  $R^{0}_{a}(D)\leq h_0(n,m,a)$ with equality if and only if $D\in \mathcal{O'}(T_{n,m})$.

\end{proof}
We obtain a  oriented  tree of matching number $m$ and order $n$ with the maximum   zeroth-order general Randi\'{c} index.
Now, we will give the maximum   zeroth-order general Randi\'{c} index for oriented unicyclic graphs with given matching number.

\begin{figure}[ht]
\begin{center}
  \includegraphics[width=14cm,height=5cm]{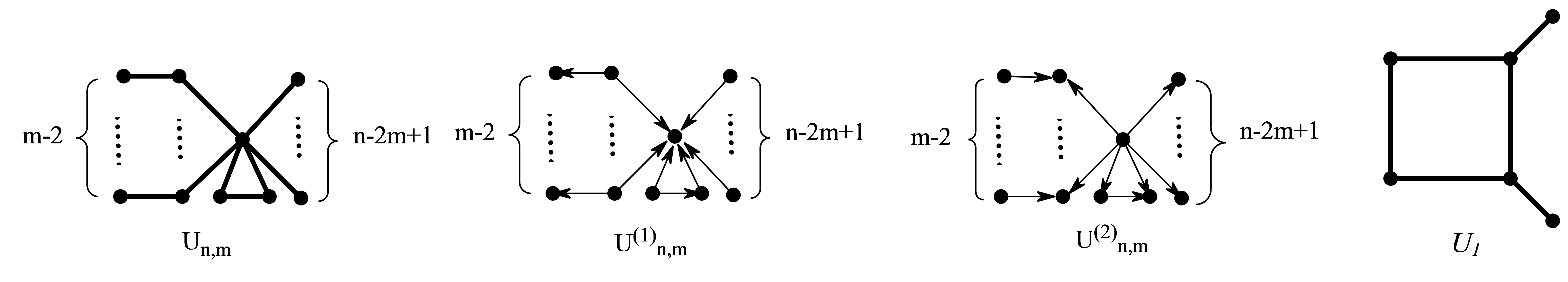}
     \end{center}
\vskip -0.5cm
\caption{{$U_1$ and Two orientations of $U_{n,m}$ :$U^{(1)}_{n,m}$,$U^{(2)}_{n,m}$.}}\label{fig-5}
\end{figure}

Let   $2\leq m\leq \lfloor\frac{n}{2}\rfloor$, $U(n,m)=\{T| T$ is a unicyclic graph  of order $n$ and matching number $m$ $\}$. Let   $U_{n,m}$, $U^{(1)}_{n,m},U^{(2)}_{n,m}$ be  shown in Figure \ref{fig-5}. Obviously, $U_{n,m}\in U(n,m)$. Let $C_n$ be the the cycle with $n$ vertices. $\tilde{U}_{n,m}=\{U^{(1)}_{n,m},U^{(2)}_{n,m}\}$.
Let

$\mathcal{U}^{*}(n,m)= \begin{cases}\tilde{\mathcal{U}}(4,2)\bigcup \{U^{(3)}(4,2), U^{(4)}(4,2)\}, & \text { if } a=1 ~and ~(n,m)=(4,2) \\ \tilde{\mathcal{U}}(6,3)\bigcup \{U^{(3)}(6,3), U^{(4)}(6,3)\}, & \text { if } a=1 ~and ~(n,m)=(6,3)\\
\tilde{\mathcal{U}}(n,m), & \text { otherwise}  \end{cases},$\\
where $U^{(3)}_{4,2},U^{(4)}_{4,2}\in \mathcal{O'}(C_4)$
and $U^{(3)}_{6,3},U^{(4)}_{6,3}\in \mathcal{O'}(U_1)$, $U_1$ is shown in Figure \ref{fig-5}.
Denote $h_1(n,m,a)=\frac{1}{2}[(-m+n+1)^{1+a}+m*2^{a+1}+n-m-2^{a+1}+1]$.
The following lemmas are  useful to prove main results.
\begin{lemma}\cite{JY223}\label{lem5}
Let $D\in \mathcal{O}(U_{4,2})$, $a\geq 1$. Then $$R^{0}_{a}(D)\leq h_1(4,2,a)$$
with equality if and only if  $D\in\{U^{(1)}_{4,2},U^{(2)}_{4,2}\}$.
\end{lemma}



\begin{figure}[ht]
\begin{center}
  \includegraphics[width=16cm,height=7cm]{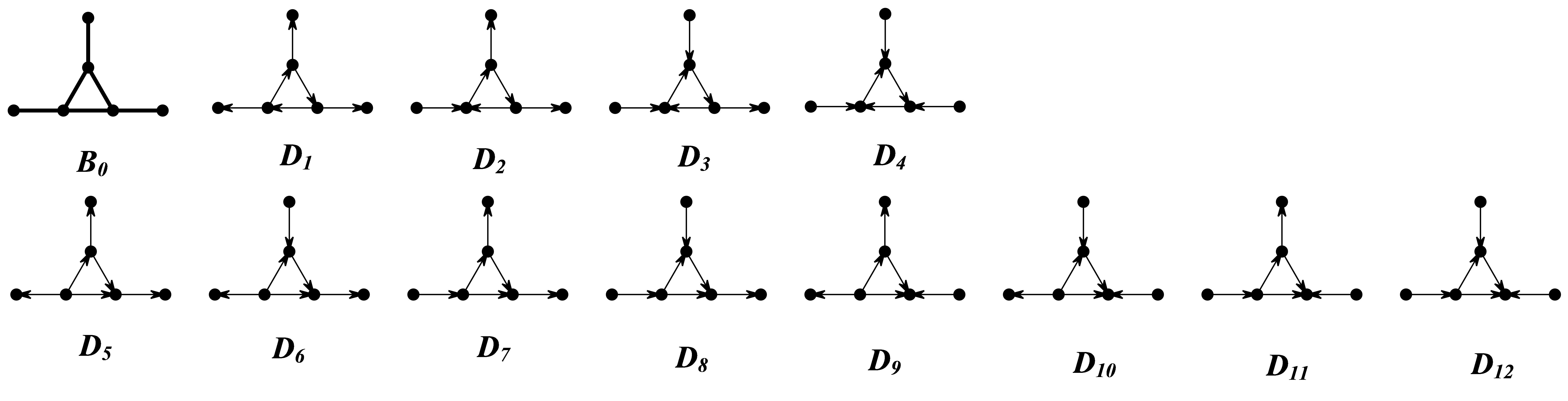}
     \end{center}
\vskip -0.5cm
\caption{ $B_0$ and $D_i$, $i=1,2,\cdot\cdot\cdot,12$ in Lemma \ref{lem7} .}\label{fig-11}
\end{figure}

\begin{lemma}\label{lem7}
  Let $B_0$ be shown in Figure \ref{fig-11}.  $D\in \mathcal{O}(B_0)$, $a\geq 1$. Then $$R^{0}_{a}(D)< h_1(6,3,a).$$
\end{lemma}

\begin{proof}
Note that $\mathcal{O}(B_0)=\{D_{i}|i=1,2,\cdot\cdot\cdot,12\}$ (see Figure \ref{fig-11}). It is easily to obtain that

$R^{0}_{a}(D_1)=R^{0}_{a}(D_2)=R^{0}_{a}(D_3)=R^{0}_{a}(D_4)=R^{0}_{a}(D_7)=R^{0}_{a}(D_8)=\frac{1}{2}(3*2^{a+1}+6),$

$R^{0}_{a}(D_5)=R^{0}_{a}(D_6)=R^{0}_{a}(D_{11})=R^{0}_{a}(D_{12})=\frac{1}{2}(2^{a+2}+3^{a+1}+5),$

$R^{0}_{a}(D_9)=R^{0}_{a}(D_{10})=\frac{1}{2}(2^{a+1}+2*3^{a+1}+4).$

Since $R^{0}_{a}(D_9)-R^{0}_{a}(D_1)=\frac{1}{2}(2*3^{a+1}-2^{a+2}-2)>0$ and
$R^{0}_{a}(D_9)-R^{0}_{a}(D_5)=\frac{1}{2}(3^{a+1}-2^{a+1}-1)>0$,
$h_1(6,3,a)-R^{0}_{a}(D_9)=\frac{1}{2}(4^{a+1}-2*3^{a+1}+2^{a+1})>0$,
 the result follows.
\end{proof}

We  consider the oriented  unicyclic graphs with a perfect matching.

\begin{theorem}\label{the8}
  Let $G \in U(2m, m)$ with $m\geq 2$, $a\geq 1$, $D\in \mathcal{O}(G)$.
 Then $$R^{0}_{a}(D)\leq h_1(2m,m,a)$$
with equality if and only if $D\in U^{*}_{2m,m} $.
\end{theorem}
\begin{proof}
 Applying induction on $m$.

   $U(4,2)=\{U_{4,2}, C_4 \}$ for $m=2$. It is easily  that $D\in \mathcal{O}(U_{4,2})$, $R^{0}_{a}(D)\leq  h_1(4,2,a)$ with equality if and only if $D\in \{U^{(1)}_{4,2},U^{(2)}_{4,2}\}$ by Lemma \ref{lem5}.  Let $f(a)=h_1(4,2,a)-2^{a+2}=\frac{1}{2}[3^{a+1}-3*2^{a+1}+3]=\frac{1}{2}[27*3^{a-2}-24*2^{a-2}+3]>0 (a\geq 2)$. As the fact that $f'(a)=\frac{1}{2}[3^{a+1}*ln3-3*2^{a+1}*ln2]>0$ $(1\leq a\leq 2)$(see Figure \ref{fig-201} of the Appendix A ) and $f(1)=0$,  $f(a)=\frac{1}{2}[3^{a+1}-3*2^{a+1}+3]\geq 0 (1\leq a\leq 2)$ with equality if and only if $a=1$. By Lemma \ref{lem1}, $D\in \mathcal{O}(C_4)$, $R^{0}_{a}(D)\leq \frac{1}{2}R^{0}_{a+1}(C_4)=2^{a+2}\leq h_1(4,2,a)$ with equality if and only if $a=1$ and $D\in \mathcal{O'}(C_4)$ which implies that $D\in \{U^{(3)}_{4,2}, U^{(4)}_{4,2}\}$.
The result follows.

 Suppose that $m \geq 3$ and the result holds for the  values smaller than  $m$. \\
\indent Let $G \in U(2m,m)$ with a perfect matching $M$.
If $G = C_{2m}$, then $D\in \mathcal{O}(C_{2m})$, by Lemma \ref{lem1} and $h_1(2m,m,a)-2^{1+a}*m=\frac{1}{2}[(m+1)^{1+a}-(1+m)*2^{a+1}+1+m]=\frac{1}{2}[(1+m)^{2}*(1+m)^{a-1}-4(1+m)*2^{a-1}+1+m]>0$, we have $R^{0}_{a}(D)\leq \frac{1}{2}R^{0}_{a+1}(C_{2m})=m*2^{a+1}< h_1(2m,m,a)$. The result follows.

If $G\neq C_{2m}$, let us consider  two cases.\\
\indent \textbf{Case 1}. $|\{v| u,v\in V(G), uv\in E(G), d_u=1,d_v=2\}|\geq 1$.

Let $d_G(u)=1$ and $ N_{G}(v)=\{w,u\}$, $d_G(w)=b$.

Since $|E(G)|-|M|=m$ and $|\{e|e\in M$ is incident with $w\}|=1$, we have $b-1\leq m$ which implies that $b\leq m+1$.
 We set $G'=G-\{v,u\}$. Then $G' \in U(2(m-1),m-1)$. Let $D'\in \mathcal{O}(G')$ such that $A(D')\bigcap A(D)=A(D')$.

By the induction hypothesis,   $R^{0}_{a}(D')\leq h_1(2(m-1),m-1,a)$. By the Lemma \ref{lem101},

\begin{equation*}
\begin{aligned}
R^{0}_{a}(D) & \leqslant R^{0}_{a}(D')+ \frac{1}{2}[1+2^{a+1}+b^{a+1}-(b-1)^{a+1}]  \\
& \leqslant h_1(2(m-1),m-1,a)+\frac{1}{2}[1+2^{a+1}+(m+1)^{a+1}-m^{a+1}] \\
& = h_1(2m,m,a)
\end{aligned}
\end{equation*}
with equality if and only if
 $D'\in U_{2(m-1),(m-1)}^{*}$ and $\{d^{-}_{D}(v)=2,d^{+}_{D}(u)=1,d_{D}^{+}(w)=m+1\}$ or $\{d^{+}_{D}(v)=2,d^{-}_{D}(u)=1, d_{D}^{-}(w)=m+1\}$, which implies that $D\in U_{2m,m}^{*}$. The result follows.

\textbf{Case 2}. $|\{v|  uv\in E(G), u,v\in V(G), d_u=1,d_v=2\}|=0$.

  $C = w_{1}w_2...w_{t}w_{1}$ denotes the unique cycle of $G$. Since $|M|=m$  , $G-V(C)$ consists of isolated vertices.

\textbf{Subcase 2.1}. For any $w_i\in V(C)$, there only exists  a pendent vertex which is adjacent to $w_i$.

 When $m\geq 8$, by Lemma \ref{lem1} and $f(a)=h_1(2m,m,a)-\frac{1}{2}[m+m*3^{a+1}]=(-1+m)*2^a+\frac{1}{2}[(1+m)^2*(1+m)^{a-1}-9m*3^{a-1}]+\frac{1}{2}>0$, we have $R^{0}_{a}(D)\leq \frac{1}{2}R^{0}_{a+1}(G)=\frac{1}{2}[m+m*3^{a+1}]<h_1(2m,m,a)$.

When $m=3$,
the result follows from  Lemma \ref{lem7}.

When $m=4$, $f(a)=3*2^{a}+\frac{1}{2}[5^{a+1}-4*3^{a+1}]+\frac{1}{2}=3*2^{a}+\frac{1}{2}[125*5^{a-2}-108*3^{a-2}]+\frac{1}{2}>0 (a\geq 2)$ and   $f(a)>0$  $(1\leq a\leq 2)$ (see Figure \ref{fig-203} of the Appendix A).   By Lemma \ref{lem1}, $R^{0}_{a}(D)\leq \frac{1}{2}[4+4*3^{a+1}]<h_1(8,4,a)$. The result follows.

When $m=5,6,7$, the proof of the cases is similar to $m=4$. The result follows.



 \textbf{Subcase 2.2}.  $|\{u|d_{G}(u)=2, u\in V(C)\}|\geq 1$.

   We suppose that $d_{w_2}=3$ and $d_{w_3}=2$. Obviously, $d_{w_1}=2$ or 3.
Let $w\in N_{G}(w_2)$ be a pendent vertex. Since $w_3$ is $M$-saturated and $w_{2}w \in M$, there exists a edge $w_{3}w_{4}\in M$, where $d_{w_4}= 2$.
 Set $T'=G-\{w_2,w\}$. Then $T'\in T(2(m-1),m-1)$. Let $D'\in \mathcal{O}(T')$ such that $A(D')\bigcap A(D)=A(D')$.

 By Lemma \ref{lem4},   we have $R^{0}_{a}(D')\leq h_0(2m-2,m-1,a)$.
 Thus
 \begin{equation*}
\begin{aligned}
R^{0}_{a}(D)
 \leqslant& R^{0}_{a}(D')+ \frac{1}{2}[(d_G(w))^{a+1}-0+(d_G(w_2))^{a+1}-0+(d_G(w_1))^{a+1}-(d_{G}(w_1)-1)^{a+1}+(d_G(w_3))^{a+1}\\
 &-(d_{G}(w_3)-1)^{a+1}]\\
 \leqslant& h_0(2m-2,m-1,a)+\frac{1}{2}[1+2*3^{1+a}-1] \\
 \leqslant& \frac{1}{2}[m-1+2^{a+1}*(m-2)+(m-1)^{1+a}+6*3^a]
\end{aligned}
\end{equation*}

 Let $g(x)=h_1(2x,x,a)-\frac{1}{2}[x-1+2^{a+1}*(x-2)+(x-1)^{1+a}+6*3^a]=-3^{a+1}+2^a+\frac{1}{2}*(1+x)^{a+1}-\frac{1}{2}*(-1+x)^{a+1}+1 (x\geq 3).$
 Since  $g'(x)=\frac{1}{2}(a+1)[(x+1)^{a}-(x-1)^{a}]>0, $
  we have $g(x)\geq g(3)$, and
 $f(a)=g(3)=-27*3^{a-2}+2^a+\frac{1}{2}(64*4^{a-2}-8*2^{a-2})+1>0 (a\geq 2).$
As the fact that $f'(a)=-3^{a+1}*ln3+2^a*ln2+\frac{1}{2}(4^{a+1}*ln4-2^{a+1}*ln2)> 0$ $(1\leq a\leq 2)$ (see Figure \ref{fig-207} of the Appendix A) and $f(1)=0,$
 $f(a)\geq 0$ with equality if and only if $a=1$.

 Consequently, $R^{0}_{a}(D) \leq h_1(2m,m,a)$ with equality if and only if $a=1$ and $D\in \{U^{(3)}_{6,3}, U^{(4)}_{6,3}\}$.
 The result follows.
\end{proof}

In order to  prove  Theorem $\ref{the11}$,  we need the following lemmas:
\begin{lemma}\cite {AY04}\label{lem9}
 Let $G\in U(n,m)$, where $n>2m$, $G \neq C_n$. Then there exist a pendent vertex $v$  and a maximum matching $M$ of $G$ such that $v$ is not  $M$-saturated.

\end{lemma}

\begin{lemma}\label{lem10}
 Let $n>2m$ and $2\leq m\leq \lfloor\frac{n}{2}\rfloor$, where $n$, $m$ are integers.

 (1)\cite{JY22} If $a=1$, then $$h_1(n,m,1)>2n.$$

 (2) If $a\geq 1$, then $$h_1(n,m,a)>2^{a}*n.$$
\end{lemma}
\begin{proof}
(2) Let $f(a)=h_1(n,m,a)-2^{a}*n=\frac{1}{2}[(1-m+n)^{a+1}+(m-1)*2^{a+1}+1+n-m]-2^{a}*n$, we have $f'(a)=\frac{1}{2}[(n-m+1)^{a+1}*ln(1-m+n)+(-n-1+m)*2^{a+1}ln2]=\frac{1}{2}[(-m+n+1)^{2}(1-m+n)^{a-1}*ln(1-m+n)+4(-1+m-n)*2^{a-1}ln2]>0$, by (1),
  $f(a)\geq f(1)>0$.
\end{proof}
Now, we consider oriented unicyclic graphs with given matching number.
\begin{theorem}\label{the11}
 Let $G\in U(n,m)$ with $2\leq m\leq \lfloor\frac{n}{2}\rfloor$, $D\in \mathcal{O}(G)$, $a\geq 1$. Then
$$R^{0}_{a}(D)\leq h_1(n,m,a)$$
with equality if and only if
$D\in U^{*}_{n,m}$.
\end{theorem}
\begin{proof}
Applying induction on $n$.

 The result follows for $n = 2m$ by Theorem \ref{the8}.

 Suppose that $n > 2m$ and the result holds for the values smaller than $n$.

  From Lemma \ref{lem1} and Lemma \ref{lem10}, if  $D\in \mathcal{O}(G)$ and $G=C_n$, then $R^{0}_{a}(D)\leq \frac{R^{0}_{a+1}(C_n)}{2}=n*2^a< h_1(n,m,a)$. The result follows.

 If $G\neq C_n$, by Lemma \ref{lem9}, there exists a pendent vertex $u$ such that  $G'=G-\{u\}\in U(n-1,m)$. Let $N_G(u)=\{v\}$, $d_G(v)=b$,  $D'\in\mathcal{O}(G')$ such that $A(D')\cap A(D)=A(D')$.

Since $|E(G)|-|M|=n-m$ and $|\{e|e\in M$ is incident with $v\}|=1$, we have $b-1\leq n-m$ which implies that $b\leq n-m+1$.

By the induction hypothesis,  $$R^{0}_{a}(D')\leq h_1(n-1,m,a).$$
By Lemma \ref{lem101}, we have
\begin{equation*}
\begin{aligned}
R^{0}_{a}(D) & \leqslant R^{0}_{a}(D')+\frac{1}{2}[1+b^{a+1}-(b-1)^{a+1}]\\
& \leqslant h_1(n-1,m,a)+ \frac{1}{2}[1+(1-m+n)^{a+1}-(-m+n)^{a+1}]\\
& = h_1(n,m,a)
\end{aligned}
\end{equation*}
with equality if and only if $R^{0}_{a}(D')= h_1(n-1,m,a)$,  $D'\in U_{n-1,m}^{*}$ and $\{d_{D}^{-}(v)=n-m+1, d_{D}^{+}(u)=1\}$ or $\{d_{D}^{+}(v)=n-m+1, d_{D}^{-}(u)=1\}$, which implies that $D\in U_{n,m}^{*}$. The result follows.
\end{proof}

\section{Maximum zeroth-order general Randi\'{c} index of orientations of bicyclic graphs with given matching number}

\begin{figure}[ht]
\begin{center}
  \includegraphics[width=14cm,height=6cm]{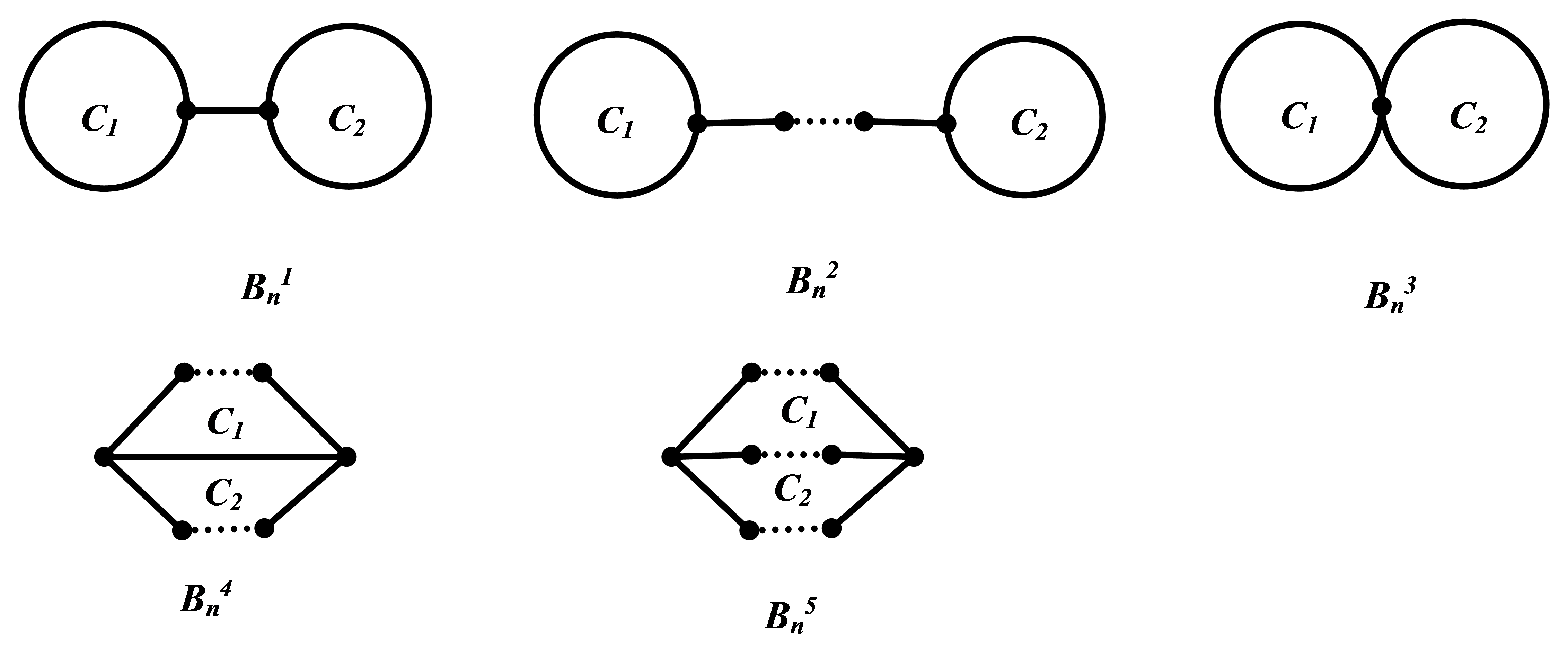}
     \end{center}
\vskip -0.5cm
\caption{ $\mathcal{B}^{i}_n, i=1,2,3,4,5$.}\label{fig-6}
\end{figure}

\begin{figure}[ht]
\begin{center}
  \includegraphics[width=10cm,height=4cm]{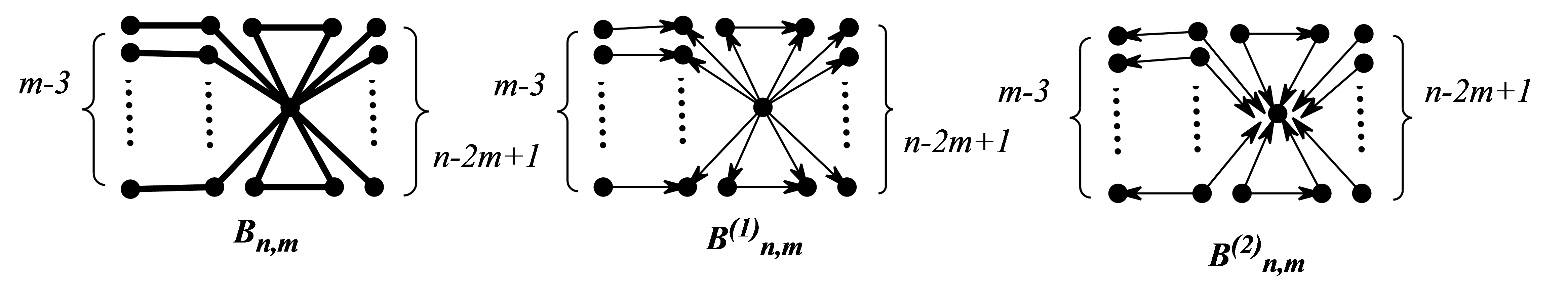}
     \end{center}
\vskip -0.5cm
\caption{ $B_{n,m}$ and its two orientations  :$B^{(1)}_{n,m}$,$B^{(2)}_{n,m}$.}\label{fig-1}
\end{figure}

In this section,
 Let $\mathcal{B}_n=\{G|G$ is a bicyclic graph  of order  $n\geq 4 \}$; $\mathcal{B}(n,m)=\{G|G$ is a bicyclic graph  of order  $n$  and matching number $m$ $\}$, where $3\leq m\leq \lfloor\frac{n}{2}\rfloor$.
 Let $\mathcal{B}^{0}_n=\{G|G$ is a bicyclic graph  of order  $n\geq 4$ and $d_{u}\neq 1$ for  any vertex $u\in V(G)\}$.
$\mathcal{B}^{0}_n$ contains five
cases according to the arrangement of cycles, which are $\mathcal{B}^{i}_n, i=1,2,3,4,5$ (see Figure \ref{fig-6}). We denote by $C_1$ and $C_2$ two independent cycles in $\mathcal{B}^{0}_n$  (see Figure \ref{fig-6}).
 It is obvious that $\mathcal{B}^{0}_n=\cup^{5}_{i=1}\mathcal{B}^{i}_{n}$.\\
\indent For $3\leq m\leq \lfloor\frac{n}{2}\rfloor$,  $B_{n,m}$, $B^{(1)}_{n,m}, B^{(2)}_{n,m}$  are shown in Figure \ref{fig-1}.
 Obviously,  $B_{n,m}\in \mathcal{B}(n,m)$.  Let $\mathcal{B}^{*}_{n,m}=\{B^{(1)}_{n,m},B^{(2)}_{n,m}\}$.

Denote $$h_2(n,m,a)=\frac{1}{2}[m*2^{a+1}-2^{a+1}-m+2+n+(n-m+2)^{a+1}].$$

\begin{lemma}\label{lem505}
Let $G$ be a connected graph, and  $w_1\in V(G)$ with $d_G(w_1)=2$, $ N_G(w_1)=\{w_2,w_3\}$,
 $d_G(w_2)\geq 2$ and $w_2w_3\notin E(G)$, $a\geq 1$. Let $G'=G-w_1w_3+w_2w_3$,
then $R^{0}_{a+1}(G')> R^{0}_{a+1}(G)$.
\end{lemma}
\begin{proof}
Let $d_G(w_i)=b_i$, $d_{G'}(w_i)=b_i'$, $i=1,2,3$.
 Obviously, $b_3=b_3'$, $b_1=b_1'+1$, $b_2=b_2'-1$, $b_2\geq 2$.
Let $f(x)=(x+1)^{a+1}-x^{a+1}$ $(x\geq 2)$, then $f'(x)=(a+1)[(x+1)^{a}-x^{a}]>0$ $(x\geq 2)$.
 Hence,
\begin{equation*}
\begin{aligned}
R^{0}_{a+1}(G)-R^{0}_{a+1}(G') &=(b_3)^{a+1}+(b_2)^{a+1}+(b_1)^{a+1}-(b_3')^{a+1}-(b_2')^{a+1}-(b_1')^{a+1}  \\
& =(b_2)^{a+1}-(b_2+1)^{a+1}+(b_1)^{a+1}-(b_1-1)^{a+1}\\
& =2^{a+1}-1-[(b_2+1)^{a+1}-(b_2)^{a+1}]<0
\end{aligned}
\end{equation*}
The result follows.
\end{proof}

\begin{figure}[ht]
\begin{center}
  \includegraphics[width=13cm,height=11cm]{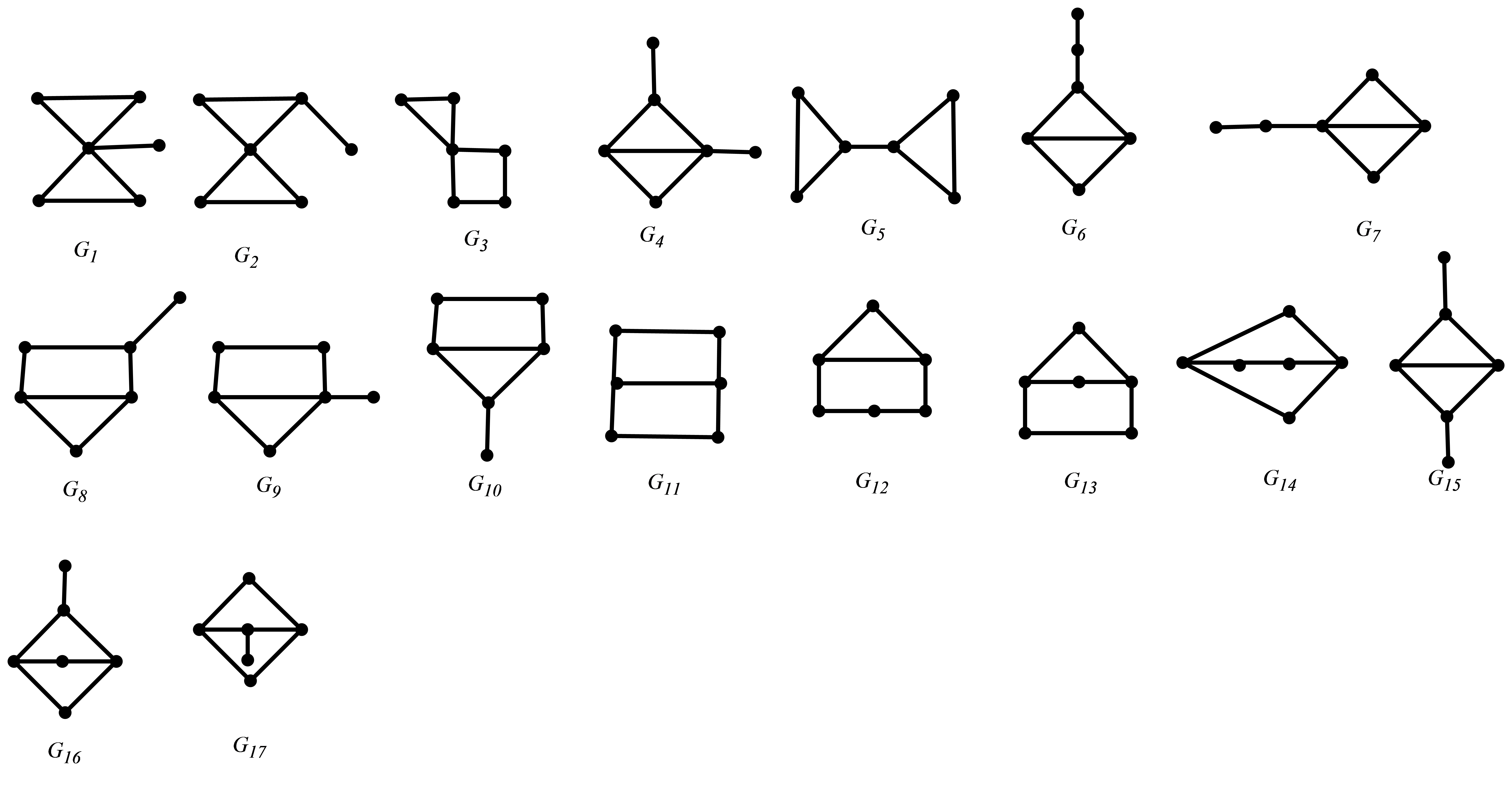}
     \end{center}
\vskip -0.5cm
\caption{{ $\mathcal{B}(6,3)$.}}\label{fig-4}
\end{figure}

\indent If $m=3$, then $\mathcal{B}(6,3)= \{G_i| i=1,2,\cdot\cdot\cdot, 17\}$ (see Figure \ref{fig-4}).
We give the maximum  zeroth-order general Randi\'{c} index for orientations of $G_1=B_{6,3}$ and $G_4$, which will be used in the proof of Theorem \ref{the502}.

\begin{lemma}\label{lem503}
Let $G_4$ be the bicyclic graph in Figure \ref{fig-4},  $D\in \mathcal{O}(G_4)$, $a\geq 1$. Then $$R^{0}_{a}(D)< h_2(6,3,a).$$

\end{lemma}

\begin{figure}[ht]
\begin{center}
  \includegraphics[width=15cm,height=8cm]{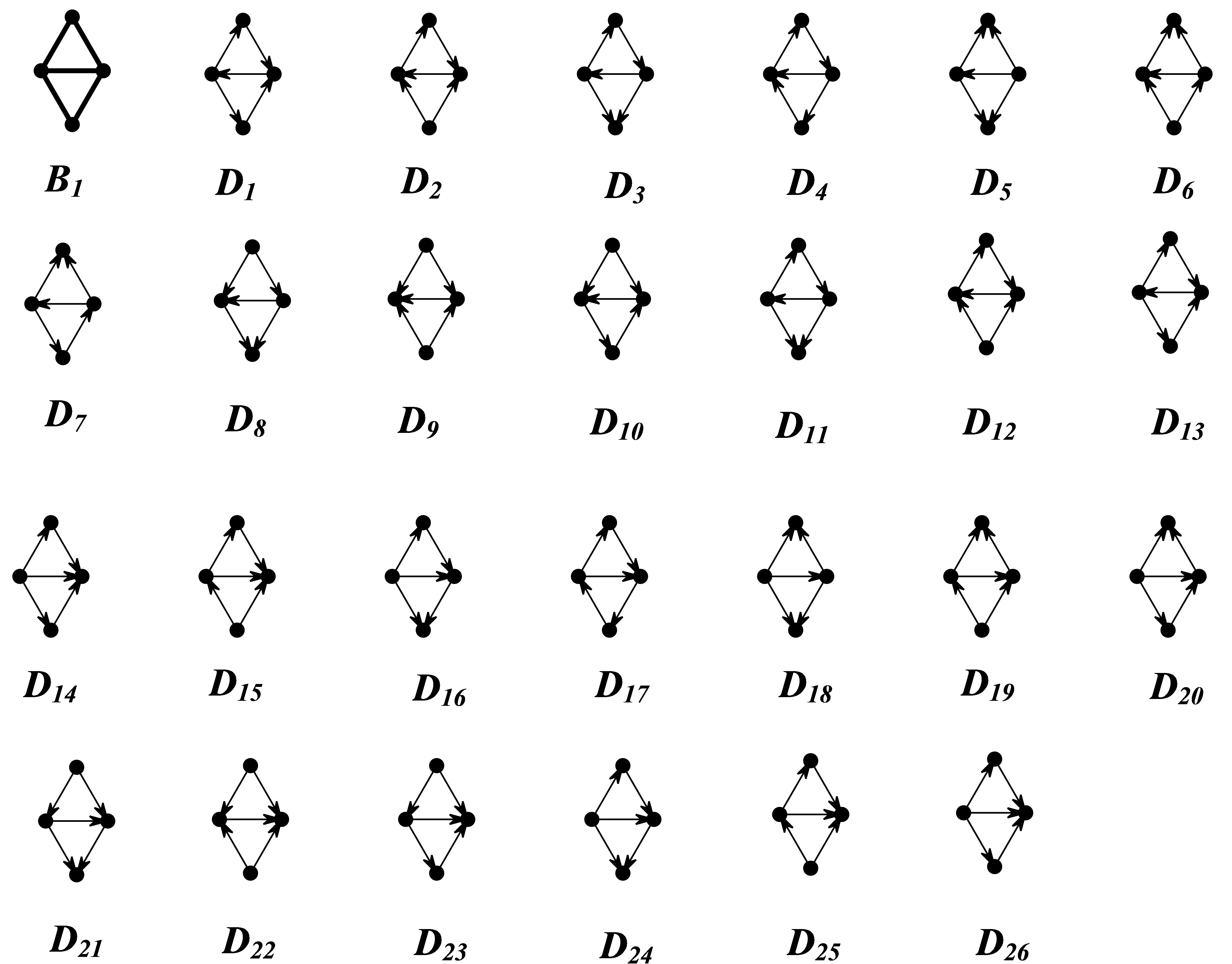}
     \end{center}
\vskip -0.5cm
\caption{{$B_1$ and its orientaitons .}}\label{fig-2}
\end{figure}

\begin{proof}
Let $B_1$ be the bicyclic graph in Figure \ref{fig-2},  $D\in \mathcal{O}(G_4)$, $D'\in \mathcal{O}(B_1)$.  $V(G_4)=\{v_4,v_3,v_2,v_1,u_2,u_1\}$,\\$V(B_1)=\{v_4,v_3,v_2,v_1\}$, and $d_{G_4}(u_i)=1$, where $i=1,2$. $u_i$ has unique neighbor $v_i$ for $i=1,2$ and $d_{G_4}(v_1)=3$,  $d_{G_4}(v_2)=4$. $d_{G_4}(v_3)=3,d_{G_4}(v_4)=2$. $d_{B_1}(v_1)=2$, $d_{B_1}(v_2)=3$. $d_{B_1}(v_3)=3,d_{B_1}(v_4)=2$.

\textbf{Case 1}. $D\in \{D| A(D)\cap A(D_1)=A(D_1), D\in \mathcal{O}(G_4)\}$.

By $1+3^{1+a}-2^{1+a}*2=9*3^{-1+a}+1-8*2^{-1+a}>0$, we have $\Sigma_{i=1}^{2}[(d^{+}_D(v_i))^{a+1}+(d^{-}_D(v_i))^{a+1}+(d^{+}_D(u_i))^{a+1}+(d^{-}_D(u_i))^{a+1}]\leq \frac{1}{2}[4+2^{a+1}+3^{a+1}]$.
Hence, $R^{0}_{a}(D)=\frac{1}{2}\Sigma_{i=1}^{2}[(d^{+}_D(v_i))^{a+1}+(d^{-}_D(v_i))^{a+1}+(d^{+}_D(u_i))^{a+1}+(d^{-}_D(u_i))^{a+1}]+\frac{1}{2}\Sigma_{i=3}^{4}[(d^{+}_D(v_i))^{a+1}+(d^{-}_D(v_i))^{a+1}] \leq \frac{1}{2}[4+3^{1+a}+2^{1+a}] +\frac{1}{2}[3+2^{1+a}]=\frac{1}{2}[7+2^{a+2}+3^{a+1}]$.

\textbf{Case 2}. $D\in \{D| A(D)\cap A(D_2)=A(D_2), D\in \mathcal{O}(G_4)\}$.

 $R^{0}_{a}(D)\leq \frac{1}{2}[5+3*2^{a+1}+3^{a+1}]$.

\textbf{ Case 3}. $D\in \{D| A(D)\cap A(D_3)=A(D_3), D\in \mathcal{O}(G_4)\}$.

 $R^{0}_{a}(D)\leq \frac{1}{2}[5+3*2^{a+1}+3^{a+1}]$.

 \textbf{Case 4}. $D\in \{D| A(D)\cap A(D_4)=A(D_4), D\in \mathcal{O}(G_4)\}$.

 $R^{0}_{a}(D)\leq \frac{1}{2}[7+2^{a+2}+3^{a+1}]$.

\textbf{ Case 5}. $D\in \{D| A(D)\cap A(D_5)=A(D_5), D\in \mathcal{O}(G_4)\}$.

 $R^{0}_{a}(D)\leq \frac{1}{2}[3+2^{a+2}+3^{a+1}+4^{a+1}]$.

 \textbf{Case 6}. $D\in \{D| A(D)\cap A(D_6)=A(D_6), D\in \mathcal{O}(G_4)\}$.

 $R^{0}_{a}(D)\leq \frac{1}{2}[4+2^{a+2}+2*3^{a+1}]$.

\textbf{ Case 7}. $D\in \{D| A(D)\cap A(D_7)=A(D_7), D\in \mathcal{O}(G_4)\}$.

 $R^{0}_{a}(D)\leq \frac{1}{2}[6+2^{a+1}+2*3^{a+1}]$.

 \textbf{Case 8}. $D\in \{D| A(D)\cap A(D_8)=A(D_8), D\in \mathcal{O}(G_4)\}$.

 $R^{0}_{a}(D)\leq \frac{1}{2}[4+2^{a+2}+2*3^{a+1}]$.

 \textbf{Case 9}. $D\in \{D| A(D)\cap A(D_9)=A(D_9), D\in \mathcal{O}(G_4)\}$.

 $R^{0}_{a}(D)\leq \frac{1}{2}[3+2^{a+1}+3^{a+2}]$.

\textbf{ Case 10}. $D\in \{D| A(D)\cap A(D_{10})=A(D_{10}), D\in \mathcal{O}(G_4)\}$.

 $R^{0}_{a}(D)\leq \frac{1}{2}[6+2^{a+1}+2*3^{a+1}]$.

 \textbf{Case 11}. $D\in \{D| A(D)\cap A(D_{11})=A(D_{11}), D\in \mathcal{O}(G_4)\}$.

 $R^{0}_{a}(D)\leq \frac{1}{2}[5+3*2^{a+1}+3^{a+1}]$.

 \textbf{Case 12}. $D\in \{D| A(D)\cap A(D_{12})=A(D_{12}), D\in \mathcal{O}(G_4)\}$.

 $R^{0}_{a}(D)\leq \frac{1}{2}[5+3*2^{a+1}+3^{a+1}]$.

 \textbf{Case 13}. $D\in \{D| A(D)\cap A(D_{13})=A(D_{13}), D\in \mathcal{O}(G_4)\}$.

 $R^{0}_{a}(D)\leq \frac{1}{2}[7+2^{a+2}+3^{a+1}]$.

Since

$\frac{1}{2}[3+2^{a+2}+4^{1+a}+3^{1+a}]-\frac{1}{2}[3+2^{a+1}+3^{a+2}]=\frac{1}{2}[2^{1+a}+4^{a+1}-2*3^{a+1}]>0$,

 $h_2(6,3,a)-\frac{1}{2}[3+2^{a+2}+4^{1+a}+3^{1+a}]=\frac{1}{2}[2+5^{a+1}-3^{1+a}-4^{a+1}]=\frac{1}{2}[2+25*5^{a-1}-9*3^{a-1}-16*4^{a-1}]>0$,

$\frac{1}{2}[3+2^{a+1}+3^{a+2}]- \frac{1}{2}[4+2^{a+2}+2*3^{a+1}]= \frac{1}{2}[-1+3^{1+a}-2^{1+a}]>0,$

$\frac{1}{2}[4+2^{a+2}+2*3^{a+1}]- \frac{1}{2}[6+2^{1+a}+2*3^{1+a}]=\frac{1}{2}[-2+2^{1+a}]>0,$

$\frac{1}{2}[6+2^{a+1}+2*3^{a+1}]-\frac{1}{2}[5+3*2^{1+a}+3^{a+1}]=\frac{1}{2}[1-2^{a+2}+3^{1+a}]=\frac{1}{2}[1-8*2^{a-1}+9*3^{a-1}]>0,$

$\frac{1}{2}[5+3*2^{1+a}+3^{1+a}]-\frac{1}{2}[7+2^{a+2}+3^{1+a}]=\frac{1}{2}[-2+2^{1+a}]>0,$\\
we have $h_2(6,3,a)>\frac{1}{2}[3+2^{a+2}+4^{1+a}+3^{1+a}]>\frac{1}{2}[3+2^{a+1}+3^{a+2}]> \frac{1}{2}[4+2^{a+2}+2*3^{a+1}]> \frac{1}{2}[6+2^{a+1}+2*3^{a+1}]>\frac{1}{2}[5+3*2^{a+1}+3^{a+1}]>\frac{1}{2}[7+2^{a+2}+3^{a+1}]$.

 It is easy to get that $D\in \{D| A(D)\cap A(D_{i})=A(D_{i}), D\in \mathcal{O}(G_4)\}$ $(i=14,15,\cdot\cdot\cdot,26)$, $R^{0}_{a}(D)<h_2(6,3,a)$.

 The result follows.
\end{proof}

\begin{lemma}\label{lem504}
Let $D\in \mathcal{O}(B_{6,3})=\mathcal{O}(G_1)$, where $G_{1}=B_{6,3}$ is shown in  Figure \ref{fig-4}, $a\geq 1$. Then $$R^{0}_{a}(D)\leq h_2(6,3,a)$$
with equality if and only if  $D\in \{B^{(1)}_{6,3},B^{(2)}_{6,3}\}$
$($see Figure \ref{fig-1}$)$.
\end{lemma}
\begin{figure}[ht]
\begin{center}
  \includegraphics[width=14cm,height=10cm]{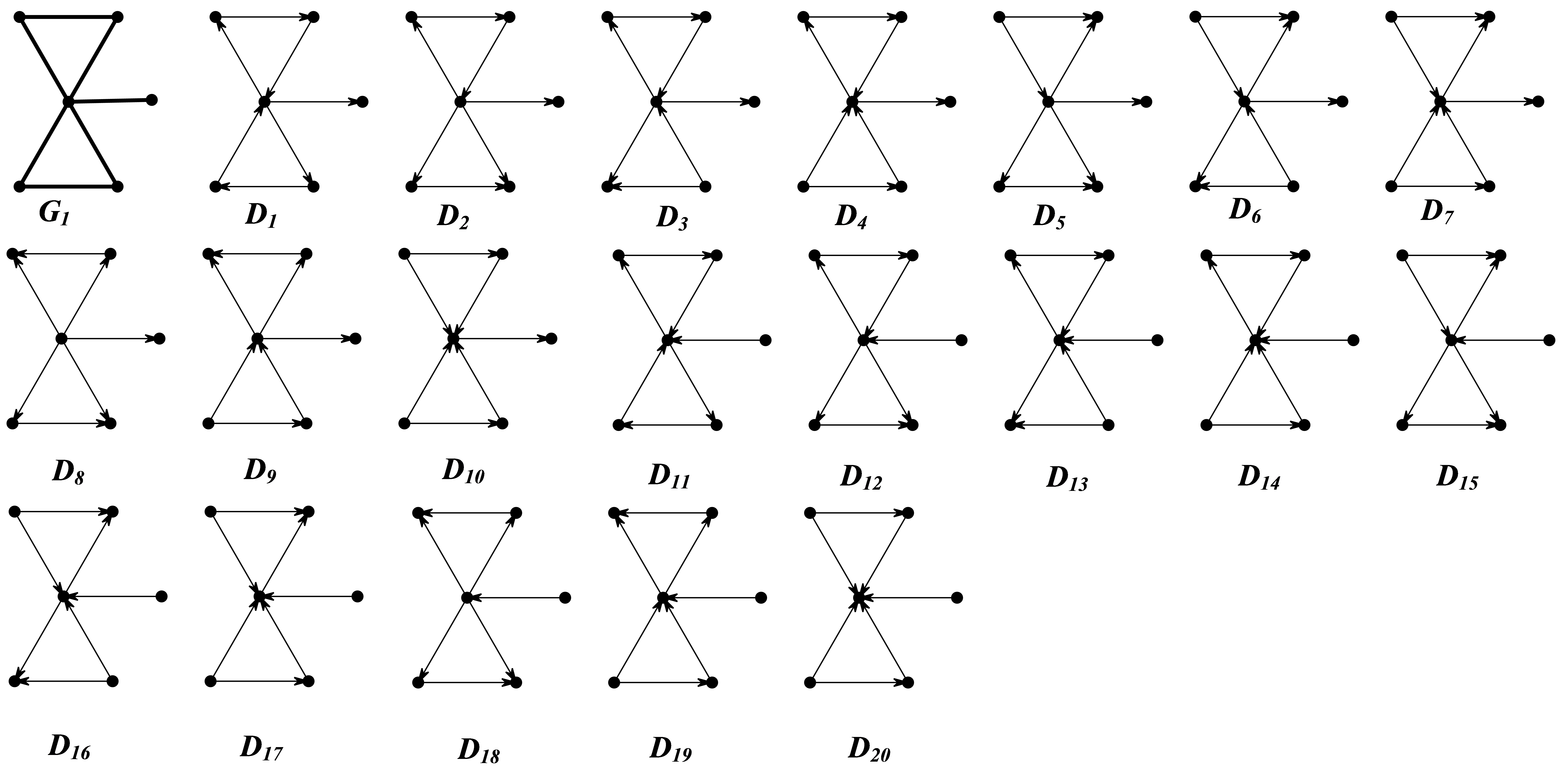}
     \end{center}
\vskip -0.5cm
\caption{{$G_1$ and its orientations.}}\label{fig-7}
\end{figure}
\begin{proof}
Let $\mathcal{O}(G_1)=\{D_1,D_2,\cdot\cdot\cdot,D_{20}\}$.
By directly calculation, we have

$R^{0}_{a}(D_1)=R^{0}_{a}(D_{11})=\frac{1}{2}(9+2^{a+1}+3^{a+1}),$

$R^{0}_{a}(D_2)=R^{0}_{a}(D_{14})=\frac{1}{2}(8+2^{a+1}+4^{a+1}),$

$R^{0}_{a}(D_3)=R^{0}_{a}(D_{9})=R^{0}_{a}(D_{13})=R^{0}_{a}(D_{19})=\frac{1}{2}(5+3*2^{a+1}+3^{a+1}),$

$R^{0}_{a}(D_4)=R^{0}_{a}(D_{12})=\frac{1}{2}(7+2^{a+2}+3^{a+1}),$

$R^{0}_{a}(D_5)=R^{0}_{a}(D_{17})=\frac{1}{2}(4+3*2^{a+1}+4^{a+1}),$

$R^{0}_{a}(D_6)=R^{0}_{a}(D_{16})=\frac{1}{2}(1+5*2^{a+1}+3^{a+1}),$

$R^{0}_{a}(D_7)=R^{0}_{a}(D_{15})=\frac{1}{2}(3+4*2^{a+1}+3^{a+1}),$

$R^{0}_{a}(D_8)=R^{0}_{a}(D_{20})=\frac{1}{2}(5+2^{a+2}+5^{a+1}),$

$R^{0}_{a}(D_{10})=R^{0}_{a}(D_{18})=\frac{1}{2}(6+2^{a+2}+4^{a+1}),$

Since $R^{0}_{a}(D_{3})- R^{0}_{a}(D_{1})=\frac{1}{2}(2^{a+2}-4)>0,$

$R^{0}_{a}(D_{3})- R^{0}_{a}(D_{4})=\frac{1}{2}(2^{a+1}-2)>0$,

$R^{0}_{a}(D_{6})- R^{0}_{a}(D_{3})=\frac{1}{2}(2^{a+2}-4)>0$,

$R^{0}_{a}(D_{6})- R^{0}_{a}(D_{7})=\frac{1}{2}(2^{a+1}-2)>0$,

$R^{0}_{a}(D_{5})- R^{0}_{a}(D_{2})=\frac{1}{2}(2^{a+2}-4)>0$,

$R^{0}_{a}(D_{5})- R^{0}_{a}(D_{10})=\frac{1}{2}(2^{a+1}-2)>0$,

$R^{0}_{a}(D_{8})- R^{0}_{a}(D_{6})=\frac{1}{2}(5^{a+1}-3^{a+1}-3*2^{a+1}+4)=\frac{1}{2}(25*5^{a-1}-9*3^{a-1}-12*2^{a-1}+4)>0$,

$R^{0}_{a}(D_{8})- R^{0}_{a}(D_{5})=\frac{1}{2}(5^{a+1}-4^{a+1}-2^{a+1}+1)=\frac{1}{2}(25*5^{a-1}-16*4^{a-1}-4*2^{a-1}+1)>0$.

The result follows.
\end{proof}
In order to get the main results, we give the maximum zeroth-order general Randi\'{c} index of orientations of bicyclic graphs with a perfect matching in which no pendent vertex has neighbor of degree 2.
\begin{lemma}\label{lem507}
Let $G\in \mathcal{B}(2m,m)$  $(m\geq 3)$, $a\geq 1$, $D\in \mathcal{O}(G)$ and no pendent vertex has neighbor of
degree 2 in $G$. Then
$$R^{0}_{a}(D)\leq h_2(2m,m,a)$$
with equality if and only if $G\in \mathcal{B}^{*}_{6,3}$.
\end{lemma}
\begin{proof}
  Suppose that $M$ is a maximum matching in $G$. Since $|M|=m$, we have that every vertex in $G$ is adjacent to at most one pendent vertex. Since $G\in \mathcal{B}(2m,m)$ and no pendent vertex has neighbor of degree 2, we  can  attach pendent vertices to  $G'\in \mathcal{B}^{0}_{k}$ ($m\leq k\leq 2m$) and get $G$. We take two cases into consideration according to that $G$ has vertices of degree 2 or not.\\
\indent  \textbf{Case 1}. $|\{v|v\in V(G), d_{G}(v)= 2\}|=0$.

Then either $k = m$ or $k =m+1$. If $k = m$, then  we can  get $G$ by attaching a pendent vertex to every vertex of  $G'\in \mathcal{B}^{0}_m$. If $k=m+1$, then  we can get $G$ by attaching a pendent vertex to every vertex of degree 2 of  $G'\in \mathcal{B}^{1}_{m+1}\cup \mathcal{B}^{4}_{m+1}\cup \mathcal{B}^{2}_{m+1}\cup \mathcal{B}^{5}_{m+1}$.

\begin{figure}[ht]
\begin{center}
  \includegraphics[width=10cm,height=6cm]{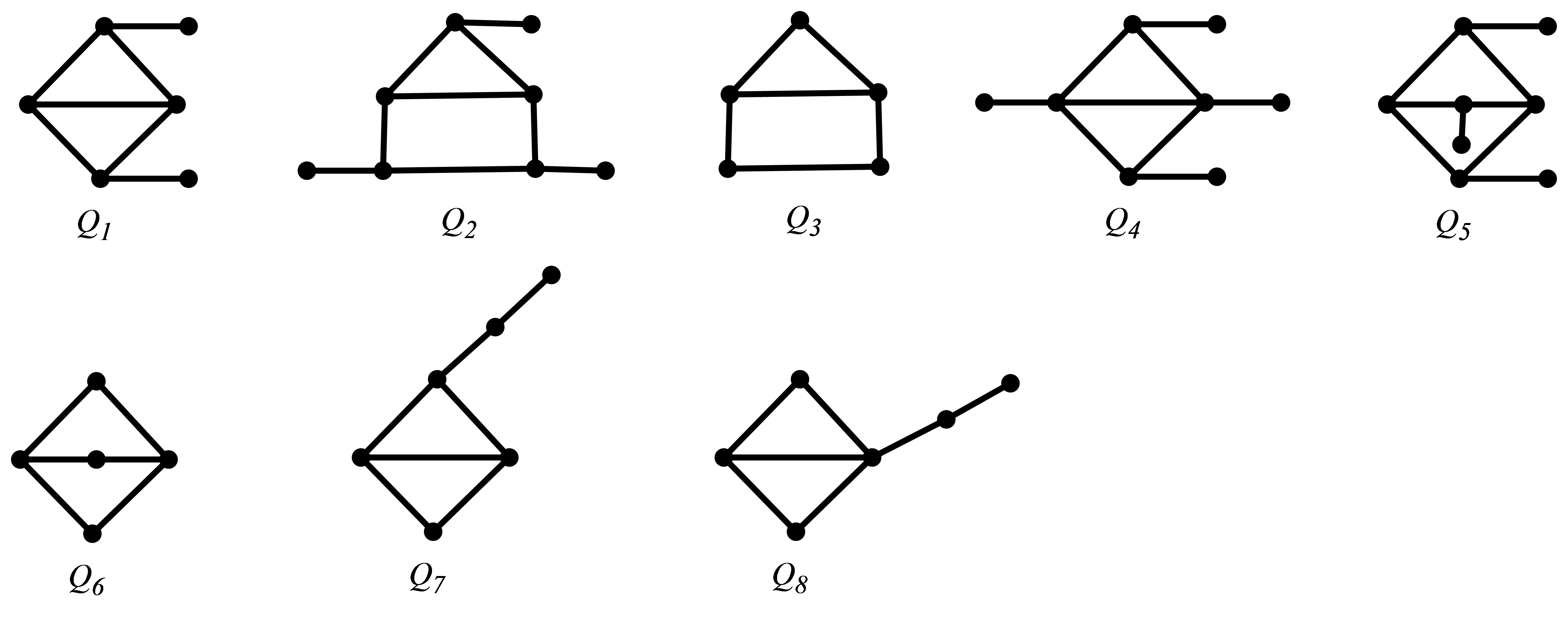}
     \end{center}
\vskip -0.5cm
\caption{{$Q_i,i=1,2,3,4,5,6,7,8$.}}\label{fig-3}
\end{figure}

\indent  If $m=3$, then  $G\cong Q_1$ (see Figure \ref{fig-3}).
Let $f(a)=h_2(6,3,a)-\frac{1}{2}[3^{a+1}*4+2]=\frac{1}{2}[2^{a+2}+5+5^{a+1}-2-12*3^a]$,
then $f(a)=\frac{1}{2}[2^{a+2}+5+125*5^{a-2}-2-108*3^{a-2}]>0 (a\geq 2)$,  $f'(a)=\frac{1}{2}[2^{a+2}*ln2+5^{a+1}*ln5-12*3^a*ln3]>0$ $(1\leq a\leq 2)$ (see Figure \ref{fig-101} of the Appendix A)  and $f(1)=0$, hence $f(a)\geq 0$ $(1\leq a\leq 2)$ with equality if and only if $a=1$.
By Lemma \ref{lem1}, $D\in \mathcal{O}(Q_1)$, $R^{0}_{a}(D)\leq \frac{1}{2}R^{0}_{a+1}(Q_1)=\frac{1}{2}[3^{a+1}*4+2]\leq h_2(6,3,a)$  with equality if and only if $a=1$ and $D\in \mathcal{O'}(Q_1)$, but by Lemma \ref{lem16} and $Q_1$ is not a bipartite graph,  $|\mathcal{O'}(Q_1)|=0$ and $R^{0}_{a}(D)< h_2(6,3,a)$. The result follows.\\
\indent  If $m=4$, we have $G'\cong Q_3$, $G\cong Q_2$ or  $G\cong Q_4$ or $G'\cong Q_6$, $G\cong Q_5$ (see Figure \ref{fig-3} ).
Let $f(a)=h_2(8,4,a)-\frac{1}{2}[3+3^{a+1}*5]=\frac{1}{2}[3*2^{a+1}+3+6^{a+1}-5*3^{a+1}],$
then $f(a)=\frac{1}{2}[3*2^{a+1}+3+216*6^{a-2}-135*3^{a-2}]> 0$ $(a\geq 2)$ and $f(a)>0$ $(1\leq a\leq 2)$ (see Figure \ref{fig-102} of the Appendix A).
By Lemma \ref{lem1}, $D\in \mathcal{O}(Q_2)$, $R^{0}_{a}(D)\leq \frac{1}{2}R^{0}_{a+1}(Q_2)=\frac{1}{2}[3+3^{a+1}*5]<h_2(8,4,a)$.

Let $f(a)=h_2(8,4,a)-\frac{1}{2}[4+2*4^{1+a}+2*3^{1+a}]=3*2^{a}+1+3*6^{a}-3^{a+1}-2^{2a+2},$
then $f(a)=3*2^{a}+1+108*6^{a-2}-27*3^{a-2}-64*4^{a-2}>0 (a\geq 2)$,  $f'(a)=3*2^{a}*ln2+3*6^{a}*ln6-3^{a+1}*ln3-2*2^{2a+2}*ln2>0$ (see Figure \ref{fig-103} of the Appendix A)  and $f(1)=0$, hence $f(a)\geq 0$ $(1\leq a\leq 2)$ with equality if and only if $a=1$. By Lemma \ref{lem1},
  $D\in \mathcal{O}(Q_4)$, $R^{0}_{a}(D)\leq \frac{1}{2}R^{0}_{a+1}(Q_4)=\frac{1}{2}[4+2*3^{a+1}+2*4^{a+1}]\leq h_2(8,4,a)$ with equality if and only if $a=1$ and $D\in \mathcal{O'}(Q_4)$, but by Lemma \ref{lem16} and $Q_4$ is not a bipartite graph,  $|\mathcal{O'}(Q_4)|=0$  and $R^{0}_{a}(D)< h_2(8,4,a)$.

Let  $f(a)=h_2(8,4,a)-\frac{1}{2}[3+5*3^{a+1}]=\frac{1}{2}[3*2^{a+1}+3+6^{a+1}-5*3^{a+1}]$, then $f(a)=\frac{1}{2}[3*2^{a+1}+3+216*6^{a-2}-135*3^{a-2}]>0 (a\geq 2)$.  $f(a)>0$ $(1\leq a\leq 2)$ (see Figure \ref{fig-666} of the Appendix A).  Hence,
$D\in \mathcal{O}(Q_5)$, $R^{0}_{a}(D)\leq \frac{1}{2}R^{0}_{a+1}(Q_5)=\frac{1}{2}[3+5*3^{a+1}]<h_2(8,4,a)$.
The result follows.\\
\indent  Now if $m\geq 5$,
then

$R^{0}_{a+1}(G)= \begin{cases}m+(m-2)*3^{a+1}+2*4^{a+1}, & \text { if } G' \in \mathcal{B}_{m}^{1} \cup \mathcal{B}_{m}^{4}\cup \mathcal{B}_{m}^{2} \cup \mathcal{B}_{m}^{5} \\ m+(m-1)*3^{a+1}+5^{a+1}, & \text { if } G' \in \mathcal{B}_{m}^{3} \\ m-1+(m+1)*3^{a+1}, & \text { if } G' \in \mathcal{B}_{m+1}^{1} \cup \mathcal{B}_{m+1}^{4}\cup \mathcal{B}_{m+1}^{2}\cup \mathcal{B}_{m+1}^{5}\end{cases}.$\\
Since $m+(m-1)*3^{a+1}+5^{a+1}-[m+(m-2)*3^{a+1}+2*4^{a+1}]=5^{a+1}-2*4^{a+1}+3^{a+1}>0$

$m+(m-1)*3^{1+a}+5^{a+1}-[m-1+(m+1)*3^{1+a}]=5^{a+1}-2*3^{a+1}+1=25*5^{a-1}-18*3^{a-1}+1>0$,\\
we have $R^{0}_{a+1}(G)\leq m+(m-1)*3^{a+1}+5^{a+1}$.

Let $f(x)=h_2(2x,x,a)-\frac{1}{2}[ x+(x-1)*3^{a+1}+5^{a+1}]=\frac{1}{2}[2+(x-1)*(2^{a+1}-3^{a+1})+(x+2)^{a+1}-5^{a+1}]$ $(x\geq 5)$,
then $f'(x)=\frac{1}{2}[2^{1+a}-3^{a+1}+(a+1)*(x+2)^{a}]> 0$,  $f(x)\geq f(5)=\frac{1}{2}[2+4(2^{a+1}-3^{a+1})+7^{a+1}-5^{a+1}]$. Since  $f(5)=\frac{1}{2}[2+4(8*2^{a-2}-27*3^{a-2})+343*7^{a-2}-125*5^{a-2}]> 0 (a\geq 2)$ and $f(5)> 0$ $(1\leq a\leq 2)$ (see Figure \ref{fig-104} of the Appendix A), we have $f(m)\geq f(5)>0$.
Thus by Lemma \ref{lem1}, $D\in \mathcal{O}(G)$, $R^{0}_{a}(D)\leq \frac{1}{2}R^{0}_{a+1}(G)\leq \frac{1}{2}[ m+(m-1)*3^{a+1}+5^{a+1}]  < h_2(2m,m,a)$. The result follows.

As the fact that
$m=3$,  $\frac{1}{2}R^{0}_{a+1}(G)\leq h_2(6,3,a)$;
$m=4$,  $\frac{1}{2}R^{0}_{a+1}(G)\leq h_2(8,4,a)$;
$m\geq 5$,  $\frac{1}{2}R^{0}_{a+1}(G)\leq h_2(2m,m,a)$,
we have $\frac{1}{2}R^{0}_{a+1}(G)\leq h_2(2m,m,a)$ in \textbf{ Case 1}.

\textbf{ Case 2}.  $|\{v|v\in V(G), d_{G}(v)= 2\}|\geq 1$.

   Let $u\in V(G)$ and $d_{G}(u)=2$, $v,w \in N_{G}(u)$ with  $d_G(v)=s\geq 2$ and $d_G(w)=
t\geq 2$. Since  $v$ and $w$ are symmetric, without lost of generality,  suppose  that $uv\in M$.\\
\indent  \textbf{SubCase 2.1}. If the cycles of $G$ have no vertex of degree 2.

Since $G$ has no
pendent vertex which has neighbor of degree 2, we have    $G'\in \mathcal{B}^{2}_{k}$
and $u$ lies on the path which is connecting two vertex-disjoint cycles of $G$. So $vw\notin E(G)$.
Since $uv\in M$ and $|M|=m$,  we have

$d_G(v)= \begin{cases}2, & \text { if $v$ lies on the path which is connecting two vertex-disjoint cycles of $G$.}   \\ s\geq 3, & \text { if $v$ lies on  cycles of $G$.} \end{cases}.$

Let
$G''=G+vw-uw$, we have $G''\in \mathcal{B}(2m,m)$. By Lemma \ref{lem505}, we have $R^{0}_{a+1}(G'')> R^{0}_{a+1}(G)$.
Obviously,

$|V'(G)|= \begin{cases}|V'(G'')|+2, & \text { if $v$ lies on the path which is connecting two vertex-disjoint cycles of $G$.}   \\ |V'(G'')|+1, & \text { if $v$ lies on  cycles of $G$.} \end{cases}$, \\
where $V'(G)$ is set of vertices of degree two in $G$ and $V'(G'')$ is set of vertices of degree two in $G''$.
 Repeating this operation from $G$ to $G''$, and we can get a bicyclic
graph  described in   \textbf{Case 1}, say $H$. By Lemma \ref{lem505} and $\frac{1}{2}R^{0}_{a+1}(G)\leq h_2(2m,m,a)$ in \textbf{ Case 1}, we have $ R^{0}_{a+1}(G)< R^{0}_{a+1}(G'')< R^{0}_{a+1}(H)\leq 2h_2(2m,m,a)$. Hence, $D\in \mathcal{O}(G)$, $R^{0}_{a}(D)<h_2(2m,m,a)$, the result follows.

  \textbf{Subcase 2.2}. We can suppose that it exists a vertex of degree 2 on some cycle of $G$, without lost of generality, say $u$.

   Let $N_G(w)=\{u=w_0, w_1,\cdot\cdot\cdot,w_{t-1}\}$, and  $G'''=G-uw$, we have $G'''\in U(2m,m)$. Note that $2\leq s,t\leq 5$ and $w$ is adjacent to at most one pendent vertex. Let  $D'\in \mathcal{O}(G''')$ such that $A(D)\bigcap A(D')=A(D')$.

By Theorem \ref{the8},   $R^{0}_{a}(D')\leq h_1(2m,m,a)$.


Let $f(x)=h_2(2x,x,a)-\frac{1}{2}[2(x+1)^{a+1}+x*2^{a+1}-x^{a+1}+x]=\frac{1}{2}[2-2^{a+1}+(x+2)^{a+1}-2(x+1)^{a+1}+x^{a+1}]$ $(x\geq 4)$, then $f'(x)=\frac{1}{2}[(a+1)(x+2)^{a}-2(a+1)(x+1)^{a}+(a+1)*x^{a}]>0$, $f(x)\geq f(4)=\frac{1}{2}[2-2^{a+1}+6^{a+1}-2*5^{a+1}+4^{a+1}]$. Let $g(a)=f(4)$, then $g(a)=\frac{1}{2}[2-2^{a+1}+1296*6^{a-3}-1250*5^{a-3}+4^{a+1}]>0 (a\geq 3)$. $g'(a)=\frac{1}{2}[-2^{a+1}*ln2+6^{a+1}*ln6-2*5^{a+1}*ln5+4^{a+1}*ln4]>0$ $(1\leq a\leq 3)$ (see Figure \ref{fig-105} of the Appendix A) and $g(1)=0$, hence $g(a)\geq 0$ $(1\leq a\leq 3)$ with equality if and only if $a=1$.
From Lemma \ref{lem101}, when $m\geq 4$, we have
\begin{equation*}
\begin{aligned}
R^{0}_{a}(D)
& \leqslant R^{0}_{a}(D')+\frac{1}{2}[(d_{G}(u))^{a+1}-(d_{G}(u)-1)^{a+1}+t^{a+1}-(t-1)^{a+1}]\\
& \leq h_1(2m,m,a)+\frac{1}{2}[2^{a+1}-1+(m+1)^{a+1}-m^{a+1}]\\
& \leq \frac{1}{2}[2(m+1)^{a+1}+m*2^{a+1}+m-m^{a+1}]\\
& \leq h_2(2m,m,a)
\end{aligned}
\end{equation*}
 with equality if and only if $a=1$ and $R^{0}_{a}(D')= h_1(2m,m,a)$, $ d_{D}^{+}(u)=2, d_{D}^{-}(w)=d_{G}(w)=5$; or $d_{D}^{-}(u)=2, d_{D}^{+}(w)=d_{G}(w)=5$.
If $D$ is satisfied with that $D'\in U^{*}_{2m,m}$, $ d_{D}^{+}(u)=2, d_{D}^{-}(w)=5$; or $d_{D}^{-}(u)=2, d_{D}^{+}(w)=5$,
then $D$ is conflicted with that no pendent vertex has neighbor of degree 2 in $G$, we have $R^{0}_{2}(D)< h_2(2m,m,a)$. The result follows.

When $m=3$,  $\mathcal{B}(6,3)/ \{G_6,G_7,G_{15}\}$ is the set of bicyclic graphs  in which  no pendent vertex has neighbor of
degree 2 and cycles contain vertex of
degree 2  (see Figure \ref{fig-4}).  Let $G\in  \mathcal{B}(6,3)/ \{G_6,G_7,G_{15}\}.$ \\
\indent  If $G=G_1$, $D\in \mathcal{O}(G_1)$, by Lemma \ref{lem504}, we have $R^{0}_{a}(D)\leq h_2(6,3,a)$
with equality if and only if  $D\in \{B^{(1)}_{6,3},B^{(2)}_{6,3}\}$.

Let $f(a)= h_2(6,3,a)-\frac{1}{2}[1+3*2^{a+1}+3^{a+1}+4^{a+1}]=\frac{1}{2}[4^{a+1}+2-2*3^{a+1}],$
then $f(a)=\frac{1}{2}[64*4^{a-2}+2-54*3^{a-2}]>0 (a\geq 2)$,  $f'(a)=\frac{1}{2}[4^{a+1}*ln4-2*3^{a+1}*ln3]>0$ $(1\leq a\leq 2)$ (see Figure \ref{fig-106} of the Appendix A)  and $f(1)=0$, hence $f(a)\geq 0$ $(1\leq a\leq 2)$ with equality if and only if $a=1$.
  If $G\in\{G_{2},G_{9}\}$ and $\frac{1}{2}R^{0}_{a+1}(G_2)=\frac{1}{2}R^{0}_{a+1}(G_9)=\frac{1}{2}[1+4^{1+a}+3^{1+a}+3*2^{1+a}]$, then by Lemma \ref{lem1}, we have  $D\in \mathcal{O}(G)$,  $R^{0}_{a}(D)\leq \frac{1}{2}R^{0}_{a+1}(G)=\frac{1}{2}[1+4^{1+a}+3^{1+a}+3*2^{1+a}]\leq h_2(6,3,a)$ with equality if and only if $a=1$ and $D\in \mathcal{O'}(G)$, but by Lemma \ref{lem16} and $G_2$, $G_9$ are not  bipartite graph,   $|\mathcal{O'}(G_2)|=0$ and $|\mathcal{O'}(G_9)|=0$, hence $R^{0}_{a}(D)< h_2(6,3,a)$. The result follows.

Let $f_{1}(a)=h_2(6,3,a)-\frac{1}{2}[4*2^{1+a}+3^{1+a}*2]=\frac{1}{2}[2*4^{1+a}-2^{a+1}+3-3^{a+2}]$, then $f_{1}(a)=\frac{1}{2}[32*4^{a-1}-4*2^{a-1}+3-27*3^{a-1}]>0.$

Let $f_{2}(a)=h_2(6,3,a)-\frac{1}{2}[1+2*2^{1+a}+3^{2+a}]=\frac{1}{2}[2*4^{1+a}+2^{a+1}+2-4*3^{1+a}],$
then $f_{2}(a)>\frac{1}{2}[128*4^{a-2}-108*3^{a-2}]> 0$ $(a\geq 2)$ and $f_{2}(a)>0$ $(1\leq a\leq 2)$ (see Figure \ref{fig-107} of the Appendix A).

Let $f_{3}(a)=h_2(6,3,a)-\frac{1}{2}[5*2^{1+a}+4^{a+1}]=\frac{1}{2}[4^{a+1}-2^{a+2}+3-3^{a+1}],$
then $f_{3}(a)=\frac{1}{2}[64*4^{a-2}-16*2^{a-2}+3-27*3^{a-2}]> 0$ $(a\geq 2)$ and $f_{3}(a)>0$ $(1\leq a\leq 2)$ (see Figure \ref{fig-108} of the Appendix A).

    If $G\in  \mathcal{B}(6,3)/ \{G_1,G_2,G_{4},G_6,G_7,G_{9},G_{15}\}$, by Lemma \ref{lem1} and $\frac{1}{2}R^{0}_{a+1}(G_5)=\frac{1}{2}R^{0}_{a+1}(G_{11})=\frac{1}{2}R^{0}_{a+1}(G_{12})=\frac{1}{2}R^{0}_{a+1}(G_{13})
=\frac{1}{2}R^{0}_{a+1}(G_{14})=\frac{1}{2}[2^{a+1}*4+2*3^{a+1}]<h_2(6,3,a)$.

$ \frac{1}{2}R^{0}_{a+1}(G_{3})=\frac{1}{2}[5*2^{a+1}+4^{a+1}]<h_2(6,3,a)$.

$\frac{1}{2}R^{0}_{a+1}(G_{8})=\frac{1}{2}R^{0}_{a+1}(G_{16})=\frac{1}{2}R^{0}_{a+1}(G_{17})
=\frac{1}{2}R^{0}_{a+1}(G_{10})=\frac{1}{2}[1+2^{a+2}+3^{a+2}]< h_2(6,3,a)$, we have  $D\in \mathcal{O}(G)$, $R^{0}_{a}(D)\leq \frac{1}{2}R^{0}_{a+1}(G)< h_2(6,3,a)$. The result follows.\\
\indent  If $G=G_4$, by Lemma \ref{lem503}, we have $D\in \mathcal{O}(G_4)$,   $R^{0}_{a}(D)< h_2(6,3,a)$. The result follows.\\
\indent  Consequently, $G\in \mathcal{B}(6,3)/ \{G_6,G_7,G_{15}\}$, $D\in \mathcal{O}(G)$,  we have $R^{0}_{a}(D)\leq h_2(6,3,a)$
with equality if and only if  $D\in \{B^{(1)}_{6,3},B^{(2)}_{6,3}\}$.
\end{proof}

Now, we are  ready to obtain the maximum zeroth-order general Randi\'{c} index of orientations of bicyclic graphs with a perfect matching.\\

\begin{theorem}\label{the502}
Let $G\in \mathcal{B}(2m,m)$, $a\geq 1$, $D\in \mathcal{O}(G)$, where $m\geq 3$. Then
$$R^{0}_{a}(D)\leq h_2(2m,m,a)$$
with equality if and only if $D\in \mathcal{B}^{*}_{2m,m}$.
\end{theorem}
\begin{proof}
Applying induction on $m$. \\
\indent  If $m=3$, then by Lemma \ref{lem507}, we
can suppose that it exists a   pendent vertex whose neighbor
degree is 2  in $G$. So we   get the bicyclic graphs $Q_7$ and $Q_8$ (see Figure \ref{fig-3}).
It is easy to get that $G\in \{Q_7, Q_8\}$, $R^{0}_{a+1}(G)\leq 1+2^{a+1}+4^{a+1}+2^{a+1}*2+3^{a+1}=1+3*2^{a+1}+3^{a+1}+4^{a+1}$.

Let $f(a)= h_2(6,3,a)-\frac{1}{2}[1+3*2^{1+a}+4^{1+a}+3^{1+a}]=\frac{1}{2}[-2^{a+1}+4+5^{1+a}-4^{a+1}-3^{1+a}],$
then $f(a)=\frac{1}{2}[-8*2^{a-2}+4+125*5^{a-2}-64*4^{a-2}-27*3^{a-2}]>0 (a\geq 2)$,  $f'(a)=\frac{1}{2}[-2^{1+a}*ln2+5^{1+a}*ln5-4^{1+a}*ln4-3^{1+a}*ln3]>0$ $(1\leq a\leq 2)$ (see Figure \ref{fig-109} of the Appendix A)  and $f(1)=0$, hence $f(a)\geq 0$ with equality if and only if $a=1$.
Then by Lemma \ref{lem1}, we have $D\in \mathcal{O}(G)$, $R^{0}_{a}(D)\leq \frac{1}{2}R^{0}_{a+1}(G)\leq \frac{1}{2}[1+3*2^{a+1}+3^{a+1}+4^{a+1}] \leq h_2(6,3,a)$  with equality if and only if $a=1$ and $D\in \mathcal{O'}(G)$, but by Lemma \ref{lem16} and $G$ is not a bipartite graph, $|\mathcal{O'}(G)|=0$ and $R^{0}_{a}(D)< h_2(6,3,a)$. The result follows.\\
 \indent  So we can suppose that $m\geq 4$ and the result follows for the vulues smaller than $m$. \\
\indent  Let $M$
be a maximum matching in $G$, we have $|M|=m$. If $G$ has no pendent vertex which has neighbor of
degree 2, then by Lemma \ref{lem507}, the result follows.\\
\indent  If $G$ has a pendent vertex $u$  whose neighbor $v$ has degree 2. Let $ N_{G}(v)=\{u,w\}$   with $d_{G}(w)=t\geq 2$, $N_{G}(w)=\{w_{0}=v,w_1,\cdot\cdot\cdot,w_{t-1}\}$,
and  $G'=G-u-v$. By $uv\in M$, we have $G'\in \mathcal{B}(2(m-1),(m-1))$.

Since $|E(G)|-|M|=m+1$ and $|\{e|e\in M$ is incident with $w\}|=1$, we have $t-1\leq m+1$ which implies that $t\leq m+2$. $D'\in \mathcal{O}(G')$ such that $A(D)\cap A(D')=A(D')$. Note that $w$ is adjacent to at most one pendent vertex in $G$.\\
\indent  By the induction hypothesis,   $R^{0}_{a}(D')\leq h_2(2m-2,m-1,a)$.

From Lemma \ref{lem101}, we have
\begin{equation*}
\begin{aligned}
R^{0}_{a}(D)&\leqslant R^{0}_{a}(D')+\frac{1}{2}[1+2^{a+1}+t^{a+1}-(t-1)^{a+1}]\\
& \leqslant h_2(2m-2,m-1,a)+\frac{1}{2}[1-(1+m)^{1+a}+(m+2)^{1+a}+2^{1+a}]\\
& =h_2(2m,m,a)
\end{aligned}
\end{equation*}
 with equality if and only if $R^{0}_{a}(D')= h_2(2m-2,m-1,a)$, $ d_{D}^{+}(v)=2, d_{D}^{-}(w)=m+2;$ $ d_{D}^{-}(v)=2, d_{D}^{+}(w)=m+2,$
 which implies that $D \in B_{2m,m}^{*}$. The result follows.
\end{proof}

 We will give the maximum zeroth-order general Randi\'{c} index of orientations of bicyclic graphs with given matching number. For this we need the following result:

\begin{lemma}\cite{YZ08}\label{lem511}
Let $G\in \mathcal{B}(n,m)$ and $G$  contains at least one pendent vertex, where $6\leq 2m< n$. Then there exists a pendent vertex $v$  and a maximum matching $M$ in $G$ such that $v$ is not $M$-saturated.
\end{lemma}

\begin{theorem}\label{the501}
Let $G\in \mathcal{B}(n,m)$, $a\geq 1$, $D\in \mathcal{O}(G)$, where $3\leq m\leq \lfloor\frac{n}{2}\rfloor$
. Then  $$R^{0}_{a}(D)\leq h_2(n,m,a)$$
with equality if and only if $D\in \mathcal{B}^{*}_{n,m}$.
\end{theorem}
\begin{proof}
Applying induction on $n$. \\
\indent  If $n=2m$, then by Theorem \ref{the502},
the result follows. \\
\indent  So we will suppose that if $n> 2m$ and the result
holds for the values  smaller than $n$. \\
\indent  If any $u\in V(G)$ and $d_{G}(u)\neq 1$, then $G\in \mathcal{B}^{0}_n$
and $n=2m+1$.

Since $2m*2^{a+1}+4^{a+1}-[(2m-1)*2^{a+1}+2*3^{a+1}]=4^{a+1}-2*3^{a+1}+2^{a+1}>0$, we have $2m*2^{1+a}+4^{a+1}>(2m-1)*2^{1+a}+2*3^{1+a}$.
Let $f(x)=h_2(n,x,a)-\frac{1}{2}[2x*2^{1+a}+4^{a+1}]=\frac{1}{2}[-(1+x)*2^{a+1}+3+x+(x+3)^{a+1}-4^{a+1}]$ $(x\geq 3)$,
then $f'(x)=\frac{1}{2}[-2^{a+1}+1+(a+1)*(x+3)^{a}]> 0$  and $f(x)\geq f(3)$.
$f(3)=\frac{1}{2}[-4*2^{a+1}+6+6^{a+1}-4^{a+1}]=\frac{1}{2}[-32*2^{a-2}+6+216*6^{a-2}-64*4^{a-2}]>0$ $(a\geq 2)$ and $f(3)>0$ $(1\leq a\leq 2)$ (see Figure \ref{fig-110} of the Appendix A).
  Since

$R^{0}_{a+1}(G)= \begin{cases}(2m-1)*2^{a+1}+2*3^{a+1}, & \text { if } $G$ \in \mathcal{B}_{2 m+1}^{1} \cup \mathcal{B}_{2 m+1}^{4} \cup \mathcal{B}_{2 m+1}^{2} \cup \mathcal{B}_{2 m+1}^{5}\\ 2m*2^{a+1}+4^{a+1}, & \text { if } $G$ \in \mathcal{B}_{2m+1}^{3}\end{cases},$\\
we have $D\in \mathcal{O}(G)$, $R^{0}_{a}(D)\leq \frac{1}{2}R^{0}_{a+1}(G) \leq \frac{1}{2}[2m*2^{a+1}+4^{a+1}]<h_2(n,m,a)$. The result follows.\\
\indent  We can suppose that $G$ contains at least one pendent vertex. Then by Lemma
\ref{lem511}, there exists a maximum matching $M$ and a pendent vertex $u\in V(G)$ such that $u$ is not $M$-saturated. Let $v\in N_{G}(u)$ and $d_{G}(v)=s\geq 2$, 
$G'=G-u$, we have $G'\in \mathcal{B}(n-1,m)$.

Since $|E(G)|-|M|=n+1-m$ and $|\{e|e\in M$ is incident with $v\}|=1$, we have $s-1\leq n+1-m$ which implies that  $s\leq n-m+2$.\\
\indent  Let  $D'\in \mathcal{O}(G')$ such that $A(D)\bigcap A(D')=A(D')$.\\
\indent  By the induction hypothesis,   $R^{0}_{a}(D')\leq h_2(n-1,m,a)$.

From Lemma \ref{lem101}, we have
\begin{equation*}
\begin{aligned}
R^{0}_{a}(D) & \leqslant R^{0}_{a}(D')+\frac{1}{2}[1+s^{a+1}-(s-1)^{a+1}]\\
& \leqslant h_2(n-1,m,a)+\frac{1}{2}[1+(n-m+2)^{a+1}-(n-m+1)^{a+1}] \\
& = h_2(n,m,a)
\end{aligned}
\end{equation*}
 with equality if and only if $R^{0}_{a}(D')= h_2(n-1,m,a)$, $\max \{d_{D}^{+}(v), d_{D}^{-}(v)\}=n-m+2$  which implies that $D\in \mathcal{B}_{n,m}^{*}$. The result follows.
\end{proof}

\textbf{Acknowledgment}.
 This work is supported by the Hunan Provincial Natural Science Foundation of China (2020JJ4423), the Department of Education of Hunan Province (19A318) and the National Natural Science Foundation of China (11971164).


Appendix A

\begin{figure}[htbp]
\centering
\begin{minipage}[t]{0.48\textwidth}
\centering
\includegraphics[width=6cm]{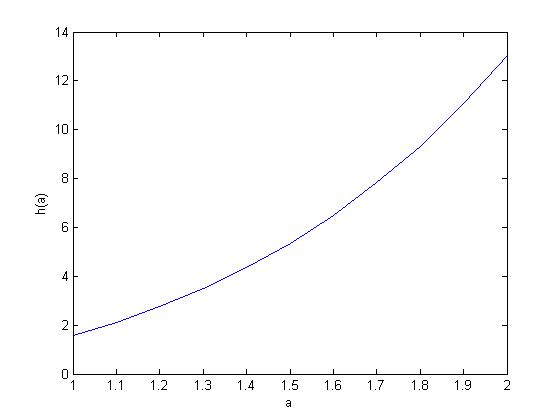}
\caption{$h(a)=2f'(a)=3^{a+1}*ln3-3*2^{a+1}*ln2$}\label{fig-201}
\end{minipage}
\begin{minipage}[t]{0.48\textwidth}
\centering
\includegraphics[width=6cm]{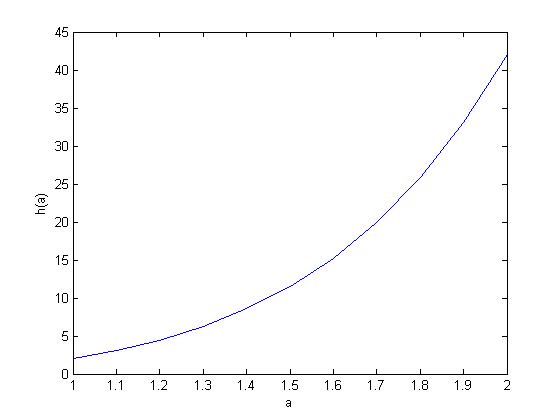}
\caption{$h(a)=2f(a)=3*2^{a+1}+5^{a+1}-4*3^{a+1}+1$}\label{fig-203}
\end{minipage}
\end{figure}

\begin{figure}[htbp]
\centering
\begin{minipage}[t]{0.48\textwidth}
\centering
\includegraphics[width=6cm]{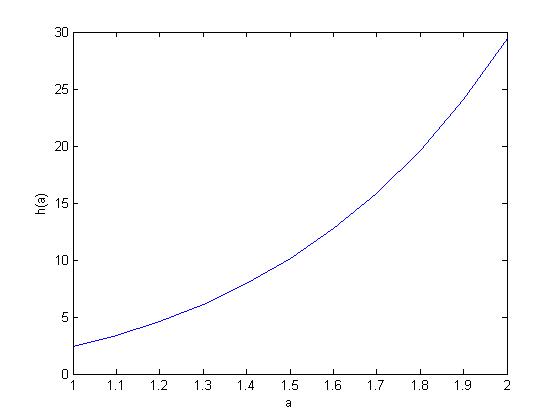}
\caption{$h(a)=2f'(a)=-2*3^{a+1}*ln(3)+4^{a+1}*ln(4)$}\label{fig-207}
\end{minipage}
\begin{minipage}[t]{0.48\textwidth}
\centering
\includegraphics[width=6cm]{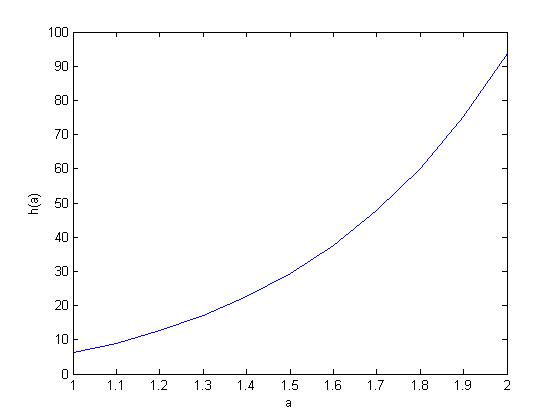}
\caption{$h(a)=2f'(a)=2^{a+2}*ln(2)+5^{a+1}*ln(5)-12*3^{a}*ln(3)$}\label{fig-101}
\end{minipage}
\end{figure}

\begin{figure}[htbp]
\centering
\begin{minipage}[t]{0.48\textwidth}
\centering
\includegraphics[width=6cm]{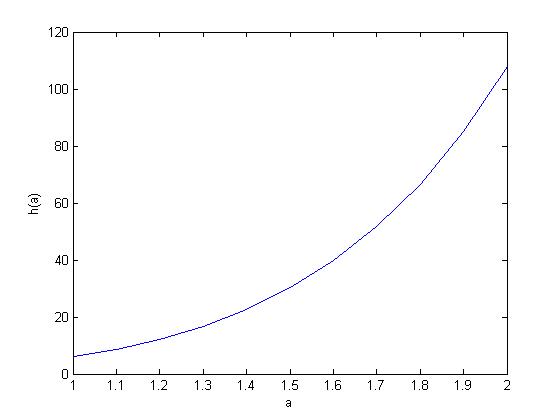}
\caption{$h(a)=2f(a)=3*2^{a+1}+3+6^{a+1}-5*3^{a+1}$}\label{fig-102}
\end{minipage}
\begin{minipage}[t]{0.48\textwidth}
\centering
\includegraphics[width=6cm]{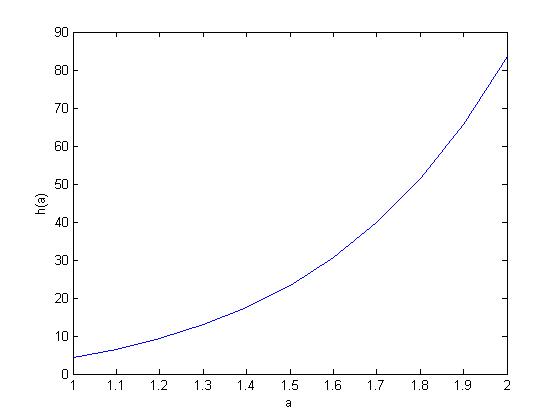}
\caption{$h(a)=f'(a)=3*2^{a}*ln(2)+3*6^{a}*ln(6)-3^{a+1}*ln(3)-2*2^{2a+2}*ln(2)$}\label{fig-103}
\end{minipage}
\end{figure}

\begin{figure}[htbp]
\centering
\begin{minipage}[t]{0.48\textwidth}
\centering
\includegraphics[width=6cm]{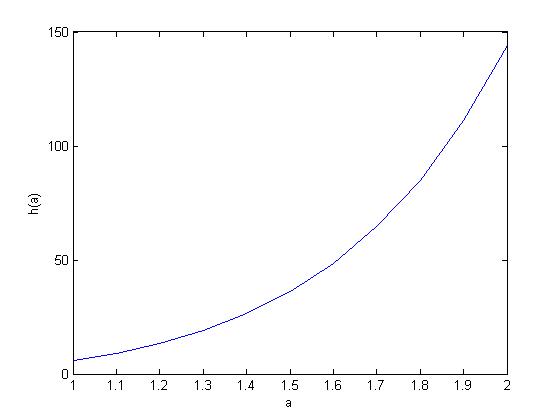}
\caption{$h(a)=2f(5)=2+4*[2^{a+1}-3^{a+1}]+7^{a+1}-5^{a+1}$}\label{fig-104}
\end{minipage}
\begin{minipage}[t]{0.48\textwidth}
\centering
\includegraphics[width=6cm]{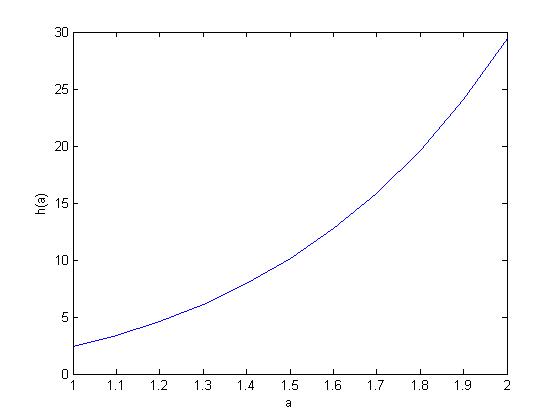}
\caption{$h(a)=2f'(a)=4^{a+1}*ln(4)-2*3^{a+1}*ln(3)$}\label{fig-106}
\end{minipage}
\end{figure}

\begin{figure}[htbp]
\centering
\begin{minipage}[t]{0.48\textwidth}
\centering
\includegraphics[width=6cm]{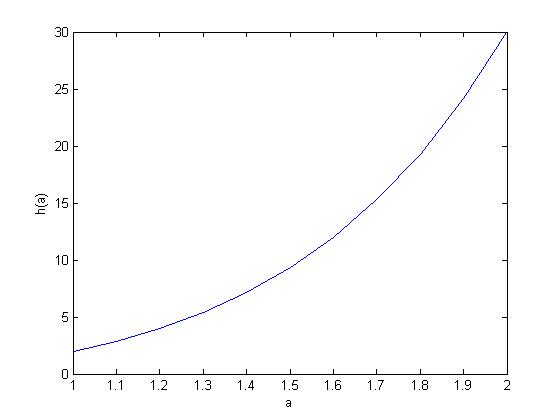}
\caption{$h(a)=2f_2(a)=2*4^{a+1}+2^{a+1}+2-4*3^{a+1}$}\label{fig-107}
\end{minipage}
\begin{minipage}[t]{0.48\textwidth}
\centering
\includegraphics[width=6cm]{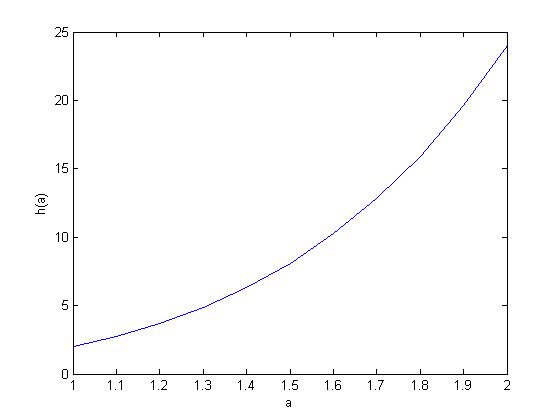}
\caption{$h(a)=2f_3(a)=4^{a+1}-2^{a+2}+3-3^{a+1}$}\label{fig-108}
\end{minipage}
\end{figure}

\begin{figure}[htbp]
\centering
\begin{minipage}[t]{0.48\textwidth}
\centering
\includegraphics[width=6cm]{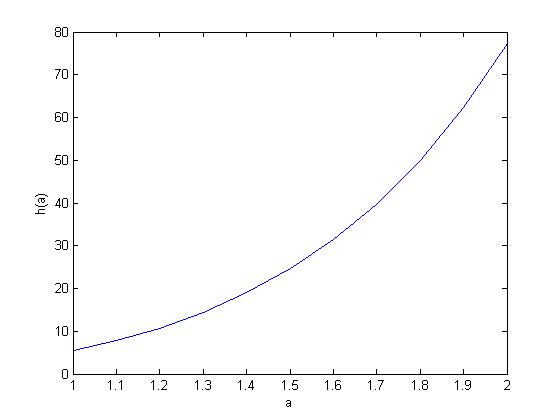}
\caption{$h(a)=2f'(a)=-2^{a+1}*ln(2)+5^{a+1}*ln(5)-4^{a+1}*ln(4)-3^{a+1}*ln(3)$}\label{fig-109}
\end{minipage}
\begin{minipage}[t]{0.48\textwidth}
\centering
\includegraphics[width=6cm]{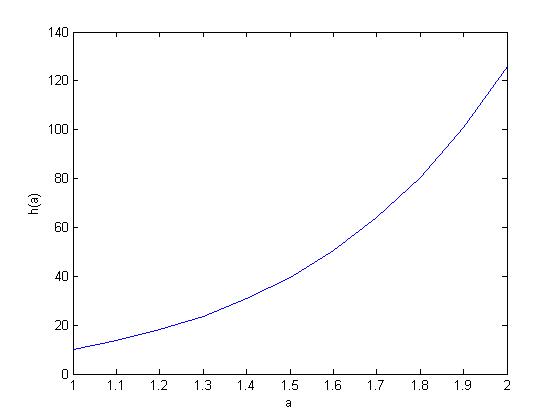}
\caption{$h(a)=2f(3)=-4*2^{a+1}+6+6^{a+1}-4^{a+1}$}\label{fig-110}
\end{minipage}
\end{figure}

\begin{figure}[htbp]
\centering
\begin{minipage}[t]{0.48\textwidth}
\centering
\includegraphics[width=6cm]{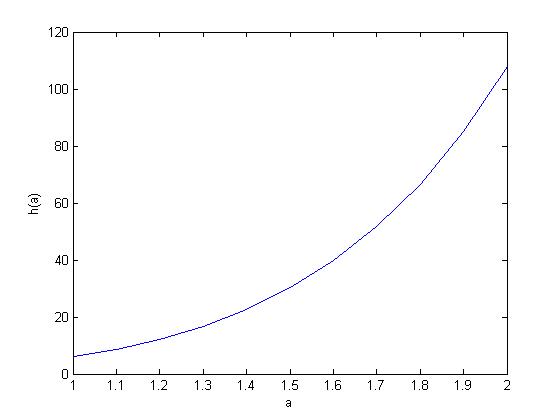}
\caption{$h(a)=2f(a)=3*2^{a+1}+3+6^{a+1}-5*3^{a+1}$}\label{fig-666}
\end{minipage}
\begin{minipage}[t]{0.48\textwidth}
\centering
\includegraphics[width=6cm]{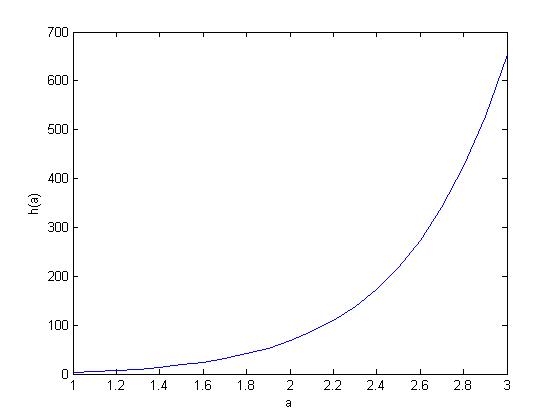}
\caption{$h(a)=2g'(a)=-2^{a+1}*ln(2)+6^{a+1}*ln(6)-2*5^{a+1}*ln(5)+4^{a+1}*ln(4)$}\label{fig-105}
\end{minipage}
\end{figure}


\begin{thebibliography}{99}
\vspace{-7pt}\bibitem{1} M. V. Diudea, QSPR/QSAR Studies by Molecular Descriptors, Nova Science Publishers, New York, NY, USA, 2001.

\vspace{-7pt}\bibitem{2} A. T. Balaban, J. Devillers, Eds., Topological Indices and Related Descriptors in QSAR and QSPAR, CRC Press,
Florida, FL, USA, 2014.

\vspace{-7pt}\bibitem{3}  X. Li, J. Zheng, A unified approach to the extremal trees for different indices, MATCH  Communications in Mathematical and in Computer Chemistry,  54  (2005) 195-208.

\vspace{-7pt}\bibitem{4} M. Liu, B. Liu, Some properties of the first general Zagreb index, Australasian Journal of Combinatorics,  47 (2010)  285.

\vspace{-7pt}\bibitem{5} J. M. Rodr\'{i}guez, J. L. S\'{a}nchez,  J. M. Sigarreta, CMMSE-on the first general Zagreb index, Journal of Mathematical
Chemistry,  56  (2018) 1849-1864.

\vspace{-7pt}\bibitem{6} H. M. Awais, M. Javaid, A. Ali, First general Zagreb index of generalized F-sum graphs, Discrete Dynamics in Nature
and Society,  2020 (2020) 1-16.

\vspace{-7pt}\bibitem{7}X. F. Pan, H. Liu, M. Liu, Sharp bounds on the zeroth-order general Randi\'{c} index of unicyclic graphs with given
 diameter, Applied Mathematics Letters,  24  (2011) 687-691.


\vspace{-7pt}\bibitem{9} M. Azari, A. Iranmanesh, Generalized Zagreb index of graphs, Studia Universitatis Babes-Bolyai
Chemia, 56  (2011) 59-70.

\vspace{-7pt}\bibitem{10} J. B. Liu, S. Javed, M. Javaid,  K. Shabbir, Computing first general Zagreb index of operations on graphs, IEEE access, 7 (2019) 47494-47502.

\vspace{-7pt}\bibitem{11} L. Bedratyuk, O. Savenko, The star sequence and the
general first Zagreb index,  MATCH  Communications in Mathematical and in Computer Chemistry, 79  (2018) 407-414.







\vspace{-7pt}\bibitem{JM212}J. Monsalve, J. Rada, Sharp
upper and lower bounds of VDB topological indices of digraphs,
Symmetry, 13  (2021) 1903.

\vspace{-7pt}\bibitem{JM211} J. Monsalve, J. Rada, Vertex-degree based topological indices of digraphs, Discrete Applied Mathematics, 295 (2021)  13-24.

\vspace{-7pt}\bibitem{JY22}J. Yang, H. Deng, Maximum first Zagreb index of orientations of unicyclic graphs with given matching number, Applied Mathematics and Computation, 427 (2022) 127131.

\vspace{-7pt}\bibitem{HD22} H. Deng, Z. Tang, J. Yang, J. Yang, M. You, On the vertex-degree based invariants of digraphs, Discrete  Mathematics Letter, 9 (2022)  2-9.

\vspace{-7pt}\bibitem{JM19} J. Monsalve, J. Rada, Oriented bipartite graphs with minimal trace norm, Linear Multilinear Algebra, 67  (2019) 1121-1131.











\vspace{-7pt}\bibitem{AY04} A. Yu, F. Tian, On the spectral radius of unicyclic graphs, MATCH  Communications in Mathematical and in Computer Chemistry,  51(2004) 97-109.







\vspace{-7pt}\bibitem{JY223}J. Yang, H. Deng, Maximum zeroth-order general Randi\'{c} index of orientations of cacti, arXiv:2205.09955 [math.CO].

\vspace{-7pt}\bibitem{YZ08} Y. Zhu, G. Liu,  J. Wang, On the Randi\'{c} index of bicyclic conjugated molecules,  in Recent Results in the Theory of Randi\'{c} Index, I. Gutman and B. Furtula, Eds. Kragujevac: University of Kragujeva, 2008, pp. 133-144.

\vspace{-7pt}\bibitem{LQ10} Q. Ling, J. Qian,  Zeroth-order general Randi\'{c} index of  trees with given matching number, Pure and Applied Mathematics, 26 (2010) 339-344.



\end{thebibliography}
\end{document}